\documentclass[a4paper,10pt]{article}
\usepackage[utf8]{inputenc}

\title{Detection of fully irreducible automorphisms in generalized Baumslag-Solitar groups}
\author{Chloé Papin}

\usepackage{graphicx}
\usepackage{amssymb}
\usepackage{tikz}
\usetikzlibrary{calc}

\usepackage{amsmath}
\usepackage[all]{xy}
\usepackage[numbers]{natbib}
\usepackage[T1]{fontenc} % Pour que les lettres accentuées soient reconnues

\usepackage{enumitem}

\usepackage{amsthm}
\newtheorem{theo}{Theorem}[section]
\newtheorem{theointro}{Theorem}
\newtheorem{prop}[theo]{Proposition}
\newtheorem{propintro}[theointro]{Proposition}

\newtheorem{lem}[theo]{Lemma}
\newtheorem{lemintro}[theointro]{Lemma}
\newtheorem{cor}[theo]{Corollary}

\theoremstyle{definition}
\newtheorem{defi}[theo]{Definition}
\theoremstyle{remark}
\newtheorem{rema}[theo]{Remark}
\newtheorem{remas}[theo]{Remarks}
\newtheorem{ex}[theo]{Examples}

\def \N{\mathbb N}
\def \Z{\mathbb Z}
\def \Q{\mathbb Q}

\def \D{\mathcal D}
\def \A{\mathcal A}
\def \E{\mathcal E}

\def \L{\mathcal L}
\DeclareMathOperator{\Aut}{Aut}
\DeclareMathOperator{\Out}{Out}

\DeclareMathOperator{\diam}{diam}
\DeclareMathOperator{\axe}{Axe}
\DeclareMathOperator{\stab}{Stab}

\DeclareMathOperator{\len}{len}
\DeclareMathOperator{\conv}{conv}
\DeclareMathOperator{\Fix}{Fix}
\DeclareMathOperator{\Isom}{Isom}

\DeclareMathOperator{\id}{id}

\DeclareMathOperator{\BS}{BS}

\DeclareMathOperator{\Lip}{Lip}

\DeclareMathOperator{\BCC}{BCC}

\DeclareMathOperator{\CV}{CV}

\DeclareMathOperator{\Wh}{Wh}

\begin{document}
  
\maketitle
 
 \newcommand{\localhost}{.}

 \begin{abstract}
  We define \emph{fully irreducible} automorphisms of generalized Baumslag-Solitar groups in analogy with fully irreducible automorphisms of free groups. We first obtain a characterization of fully irreducible automorphisms analogous to a condition given by Kapovich. Next we discuss the existence of \emph{pseudo-periodic} conjugacy classes based on the study of Nielsen paths in a train track representative. Finally we give an algorithm which finds out whether an automorphism with a train track representative is pseudo-atoroidal and fully irreducible.
 \end{abstract}

\section*{Introduction}

Baumslag-Solitar groups are defined by
\[
 \BS(p,q) = \langle a, t | t a^p t^{-1} = a^q \rangle
\]
They were introduced by Baumslag and Solitar in \cite{BaumslagSolitarSomeTwoGenerator} as examples of non-Hopfian groups.

\begin{figure}
 \centering
 \includegraphics[scale=0.4]{\localhost/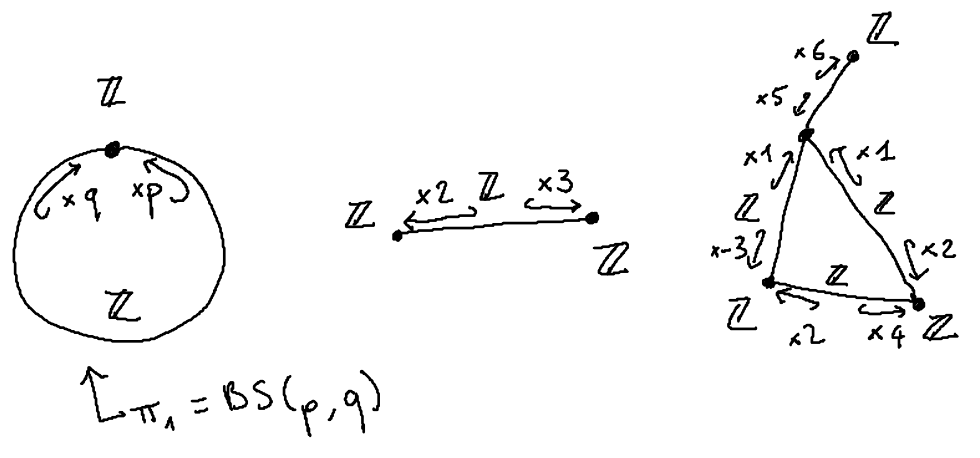}
 \caption{Examples of graphs of groups with cyclic vertex and edge groups} \label{fig:gdgBS}
\end{figure}

Generalized Baumslag-Solitar (GBS) groups are defined as fundamental groups of finite graphs of infinite cyclic groups (see Figure \ref{fig:gdgBS} for examples). Equivalently they are groups which act on a simplicial tree with infinite cyclic edge and vertex stabilizers and finite quotient.

In this paper we are interested in the outer automorphism group $\Out(G)$ for a GBS group $G$.

\bigskip

GBS groups naturally act on trees with infinite cyclic edge and vertex stabilizers. In general a GBS group admits infinitely many such actions, and there are infinitely many corresponding graphs of groups. In \cite{ForesterSplittings} Forester introduced \emph{deformation spaces} associated to a GBS group $G$. It consists in the space of $G$-trees which have the same elliptic groups as a given tree $T$. Here we consider the \emph{cyclic deformation space} which is the space of all minimal actions of $G$ on simplicial trees with infinite cyclic edge and vertex stabilizers, up to $G$-equivariant isomorphism and homothety (unless $G$ is isomorphic to $\Z$, $\Z^2$ or the fundamental group of a Klein bottle). The cyclic deformation space $\D$ is the GBS analogue of the outer space of Culler and Vogtmann $\CV_N$ for the free group $F_N$: outer space is a contractible space with a proper action of $\Out(F_N)$.

For a GBS group $G$, there is a natural action of $\Aut(G)$ on $\D$ by pre-composition of the action of $G$: if $\phi \in \\Aut(G)$ and $(T, \rho) \in \D$ where $\rho: G \to \Isom(T)$ is the action, then $(T, \rho) \cdot \phi = (T, \rho \circ \phi)$. Since inner automorphisms act trivially this is actually an action of $\Out(G)$. In the rest of the paper, we distinguish actual automorphisms in $\Aut(G)$ and their outer class in $\Out(G)$: we denote the former with lower case letters and the latter with upper case letters.

Studying actions of $G$ on trees with cyclic stabilizers is equivalent to studying \emph{marked graphs of groups} $(\Gamma, \sigma)$ where $\Gamma$ is a graph of cyclic groups and $\sigma$ is an isomorphism $G \simeq \pi_1(\Gamma)$.

An automorphism $\Phi \in \Out(F_N)$ is called \emph{fully irreducible} when no conjugacy class of free factors is $\Phi$-periodic. In the case of GBS groups, free factors are not relevant but there is an analogue called \emph{special factor}, which we develop in \cite{PapinWhitehead}. 
A special factor is a subgroup of $G$ which is the fundamental group of a subgraph of groups in some graph of cyclic groups for $G$. We can define a fully irreducible automorphism $\Phi \in \Out(G)$: it is an automorphism which has no periodic conjugacy class of special factor.

In \cite{BestvinaHandelTrainTracks} Bestvina and Handel prove that any fully irreducible element $\Phi \in \Out(F_N)$ admits a \emph{train track representative}, which is an action on a tree $T\in\D$ with a $\phi$-equivariant map $f: T \to T$ for some $\phi \in \Phi$ whose iterates stretch equally all edges and send them on geodesic paths. The set of train track trees behave like a translation axis for $\Phi$ in $\CV_N$. Train track representatives are a key tool for studying fully irreducible automorphisms and they can be computed. However it is not known whether train track representatives always exist for fully irreducible automorphisms of a GBS group in general, although there are many examples. When $G$ has no \emph{non-trivial integer modulus} (see for example \cite{LevittAutomorphisms} for a definition of modulus) then there is a bound on the number of edges in graphs of cyclic groups for $G$ (\cite{ForesterSplittings}). The results in \cite{MeinertTheLipschitzMetric} imply that fully irreducible automorphisms of $G$ in this case admit train track representatives. Yet we do not know whether such groups always admit interesting fully irreducible automorphisms. 

Another example is $G =\BS(p,pn)$ for $p > 1$. In \cite{BouetteThese} Bouette proved that in that case, all automorphisms are reducible since there is a globally invariant conjugacy class of special factors. However by replacing $\D$ by another deformation space allowing greater vertex stabilizers, namely those in the invariant conjugacy class, there are fully irreducible automorphisms and they admit train tracks.

In this paper we adapt an algorithm given by Kapovich in \cite{KapovichAlgorithmic} and \cite{KapovichDetectingFully} which, given an automorphism $\Phi$ which is \emph{atoroidal}, i.e. has no periodic conjugacy class, and a train track representative for $\Phi$, decides whether or not $\Phi$ is fully irreducible. The atoroidal condition can actually be weakened. In the original papers by Kapovich the algorithm applies to any automorphism since there is an algorithm in \cite{BestvinaHandelTrainTracks} which either finds a train track representative for $\Phi$ or gives a proof that $\Phi$ is not reducible.

A non-solvable GBS group is a GBS group which is neither isomorphic to $\Z$, $\Z^2$, the fundamental group of a Klein bottle nor $\BS(1,n)$ for some $n \in \Z \setminus \{0\}$.

A conjugacy class $[g] \in G$ is \emph{loxodromic} if there exists $T \in \D$ (equivalently for any $ T \in \D$) the action of $g$ on $T$ is loxodromic.

A loxodromic conjugacy class $[g]$ with $g \in G$ is \emph{pseudo-periodic} for $\Phi \in \Out(G)$ if $\|\phi^n(g)\|$ is bounded for $n \in \N$, for some automorphism $\phi \in \Phi$. While in a free group context this would ensure that the conjugacy class of $g$ is actually periodic, it is not the case here. An outer class $\Phi \in \Out(G)$ is \emph{pseudo-atoroidal} if it has no pseudo-periodic hyperbolic conjugacy class.

 The paper is centred around the two following theorems, which are independant. 
\begin{theointro} \label{theo:mainTheo}
 There is an algorithm which takes a non-solvable GBS group $G$, a pseudo-atoroidal automorphism $\Phi \in \Out(G)$ and a train track representative $f:T \to T$ for $\Phi$, and decides whether $\Phi$ is fully irreducible.
\end{theointro}
\begin{theointro} \label{theo:decider-patoroidal}
 There is an algorithm which takes a non-solvable GBS group $G$ and an outer automorphism $\Phi \in \Out(G)$ with a train track representative $f : T \to T$ and decides whether $\Phi$ is pseudo-atoroidal.
\end{theointro}

\paragraph{Structure of pseudo-periodic conjugacy classes and subgroups.}

The fixed group for an automorphism $\phi \in \Aut(G)$ is the subgroup of $G$ whose elements are fixed for $\phi$. When $\phi$ is an automorphism of a free group, the fixed group of $\phi$ is always finitely generated (\cite{GerstenFixedPoints}).
In GBS groups, fixed groups are not always finitely generated: see Example \ref{ex:groupe-fixe}. The notion of pseudo-periodic subgroup associated to an automorphism $\phi \in \Aut(G)$ seems to be more interesting.

\begin{ex}\label{ex:groupe-fixe}
 Let $G:= \BS(2,4) = \langle a, t| ta^2 t^{-1} = a^4 \rangle$.
 Let $\phi \in \Aut(G)$ be the inner automorphism defined by $g \mapsto a^2 g a^{-2}$. The fixed group $\Fix(\phi)$ is the subgroup generated by $\{ t^{-k} a^2 t^k, k \in \Z\}$ and it is not finitely generated. In fact, let $T \in \D$ be the standard tree for $\BS(2,4)$ of Figure \ref{fig:gdgBS}. Let $v$ be the unique vertex of $T$ with stabilizer $\langle a \rangle$. Let $e$ be the edge from $v$ to $tv$. 
 The minimal invariant subtree for $\Fix(\phi)$ is 
 the subset of $T$ which can be reached from $v$ by an edge path $\gamma:= e_1, \dots, e_n$ such that for every $i \in \{1, \dots, n\}$ the subpath $e_1, \dots, e_i$ does not contain more edges in $G \cdot e$  than in $G \cdot \bar e$.
 
 The quotient of the minimal subtree by $\Fix( \phi)$ has infinite diameter. Therefore $\Fix(\phi)$ is not finitely generated.
\end{ex}

In free groups, pseudo-periodic conjugacy classes are actually periodic. However, in a GBS group, this is not automatically true. In \cite{PapinWhitehead} we proved the existence of a minimal special factor containing $g$ for every loxodromic $g \in G$. A conjugacy class is \emph{simple} if its minimal special factor is not $G$. 
If $[g]$ is pseudo-periodic, then the conjugacy class of the minimal factor containing $g$ is periodic: if $g$ is simple, then $G$ has a proper periodic conjugacy class of special factors. Thus we have:
 \begin{lemintro}
  If there exists a simple pseudo-periodic conjugacy class for $\Phi \in \Out(G)$, then $\Phi$ is not fully irreducible.
 \end{lemintro}

 The strategy for finding pseudo-periodic conjugacy classes is based on \emph{Nielsen paths}, just as the study of periodic conjugacy classes in free groups. Let $f: T \to T$ be a train track representative for $\Phi \in \Out(G)$. A \emph{Nielsen path} is a finite path $[x,y]$ in $T$ such that there exists $g \in G$ such that $f(x)= gx$ and $f(y)=gy$. More generally a periodic Nielsen path is a Nielsen path for $f^n$ for some $n \geq 1$. They can be split into concatenations of \emph{periodic indivisible Nielsen paths} (pINP). The link between pseudo-periodic conjugacy classes for $\Phi$ and periodic indivisible Nielsen paths for $f$ is given by:
 \begin{propintro} \label{prop:concatenationNielsenIntro}
  Let $g \in G$ be an element whose conjugacy class is pseudo-periodic for $\Phi \in \Out(G)$. Suppose $f: T \to T$ is a train-track representative for $\Phi$. Then the axis of $g$ in $T$ is a concatenation of periodic indivisible Nielsen paths.
 \end{propintro}
 
 Consider a periodic indivisible Nielsen path $\eta \subset T$. Define $VY(\eta)$ as the set of points of $T$ which can be joined to one endpoint of $\eta$ by a concatenation of periodic indivisible Nielsen paths. By Proposition \ref{prop:concatenationNielsenIntro} the axis of $g \in G$ such that the conjugacy class $[g]$ is pseudo-periodic is contained in the convex hull of $VY(\eta)$ for some periodic indivisible Nielsen path $\eta$ and that that $g \in \stab(VY(\eta))$.
 
 Like in the case of free groups, there are finitely many $G$-orbits of periodic indivisible Nielsen paths, and they can be computed:
 \begin{theointro} \label{theo:compute-pinps-intro}
 Let $\Phi \in \Out(G)$ and let $f: T \to T$ be a train track map for $\Phi$. There are finitely many orbits of periodic indivisible Nielsen paths for $f$. 
 
 Furthermore there is an algorithm which finds all orbits of pINPs.
\end{theointro}
 Consequently there exist finitely many sets $VY(\eta)$ where $\eta$ is a pINP up to translation, so there exist finitely many conjugacy classes of subgroups $\stab(VY(\eta))$ which could contain elements whose conjugacy class is pseudo-periodic:
 \begin{theointro}
  Let $\Phi \in \Out(G)$. There exist finitely many subgroups 
  \[
  G_1=\stab(VY(\eta_1)), \dots, G_k= \stab(VY(\eta_k))
  \]
  in $G$, well defined up to conjugacy, for which one can compute a finite set of generators, such that:
  \begin{itemize}
   \item If a conjugacy class $[g]$ with $g \in G$ is pseudo-periodic for $\Phi$ then there exists $i \in \{1, \dots, k\}$ such that a conjugate of $g$ belongs to $G_i$.
   \item Conversely, for all $i \in \{1, \dots, k\}$, all loxodromic elements of $G_i$ have a pseudo-periodic conjugacy class.
  \end{itemize}
 \end{theointro}
 
 Theorem \ref{theo:decider-patoroidal} follows.
 
 \bigskip
 
 Given an actual automorphism $\phi \in \Aut(G)$ and a tree $T \in \D$, one can define the \emph{pseudo-periodic subgroup} associated to $G_\phi$ as the subgroup
 \[
  G_\phi:=\{g \in G/  d_T(*, \phi^n(g) *) \text{ is bounded for } n \in \N \}
 \]
 where $*$ is any base point in $T$. The group $G_{\phi}$ does not depend on the choice of the base point $*$, nor on the choice of $T \in \D$.

Pseudo-periodic subgroups of GBS groups are related to fixed groups for automorphisms of free groups. There is also a link with Nielsen paths and pseudo-periodic conjugacy classes:
\begin{theointro} 
 Let $\Phi \in \Out(G)$ and let $f: T \to T$ be a train track representative for $\Phi$. Let $\psi \in \Aut(G)$ such that $\psi \in \Phi^k$ for some $k \in \N$. If $G_{\psi}$ contains a loxodromic element, there exists a pINP $\eta \subset T$ such that
 \[
  G_{\psi} \subset \stab(VY(\eta))
 \]
 
 Moreover, for any pINP $\eta$ such that $\stab(VY(\eta))$ contains a loxodromic element, there exists $k \in \N$ and $\psi \in \Phi^k$ such that $\stab(VY(\eta)) = G_{\psi}$.
\end{theointro}
This indicates that the subgroups $\stab(VY(\eta))$ for $\eta$ a pINP are the maximal subgroups among the pseudo-periodic subgroups which contain loxodromic elements.

\paragraph{Periodic special factors and stable lamination.}

Given an outer automorphism $\Phi \in \Out(G)$ with a train track representative $f : T \to T$, one can define the \emph{stable lamination} $\Lambda^+$ associated to $T$ as the $G$-invariant collection of bi-infinite geodesics of $T$ obtained by taking the closure of the set of translates of iterates $f^n(e)$ for an edge $e \in T$.

Theorem \ref{theo:mainTheo} is a consequence of the following criterion extending \cite{KapovichAlgorithmic}: 
\begin{theointro} \label{theo:critere-iwip}
 Let $G$ be a non-solvable generalized Baumslag-Solitar group. 
 Let $\Phi \in \Out(G)$ be an automorphism of $G$ with a primitive train track $f: T \to T$ with no simple pseudo-periodic element. Then $\Phi$ is fully irreducible if and only if all Whitehead graphs $\Wh_T(\Lambda^+,v)$ are connected.
\end{theointro}
 A train track map is \emph{primitive} when the associated transition matrix is primitive.

 We actually prove a slightly stronger version of this criterion, using restricted deformation spaces instead of the standard deformation space (see \cite[Def. 3.12]{GuirardelLevitt07}): restricted deformation spaces have a restriction on allowed edge groups as well as on vertex groups. 
 
 From this criterion we deduce the algorithm of Theorem \ref{theo:mainTheo}.
 The algorithmicity here comes from the fact that $\Lambda^+(f)$ is quasi-periodic, which guarantees that the computation of the turns taken by the lamination takes only finitely many steps. The computation of Whitehead graphs is possible since the trees in $\D$ are locally finite.

 \bigskip
 
 Note that the criterion given by Theorem \ref{theo:critere-iwip} does not require that $\Phi$ be pseudo-atoroidal: Theorem \ref{theo:mainTheo} actually applies to automorphisms with non-simple pseudo-periodic conjugacy classes. However we do not know how to algorithmically determine whether an automorphism satisfies this hypothesis. We can decide whether there exist pseudo-periodic conjugacy classes, and for a given element we can decide whether it is algorithmic using the results of \cite{PapinWhitehead} but the set of pseudo-periodic conjugacy classes could be complicated.
 
 In order to answer this question, we would need to know whether one of the subgroups $\stab(VY(\eta))$ defined before contains a simple element. Without any algorithmic answer to this question, we need the pseudo-atoroidal condition to assemble Theorems \ref{theo:decider-patoroidal} and \ref{theo:mainTheo}.
 
\bigskip

In Section \ref{sec:iwips-definitions} we define irreducible automorphisms and train track representatives, as well as the stable lamination which comes with a train track map. In Section \ref{sec:detecter} we prove Theorem \ref{theo:critere-iwip}. Section \ref{sec:nielsen} is about pseudo-periodic conjugacy classes and how they may be understood using Nielsen paths; pseudo-periodic subgroups and their computation is developed in Section \ref{sec:groupes-pseudo}. These two sections can be read independently of Section \ref{sec:detecter}. Finally Section \ref{sec:conclu} puts together the arguments for Theorem \ref{theo:mainTheo}.

\section{Special factors, automorphisms of GBS groups} \label{sec:iwips-definitions}

A Generalized Baumslag-Solitar group (GBS) is a group which admits a minimal action on a simplicial tree with infinite cyclic vertex and edge stabilizers. Equivalently it is a group isomorphic to the fundamental group of a finite graph of infinite cyclic groups on both vertices and edges.

GBS groups isomorphic to $\Z$, $\Z^2$ and the fundamental group of a Klein bottle are called \emph{elementary}. Along with Baumslag-Solitar groups $\BS(1,n)$ for $n \in \Z \setminus \{0\}$, they are the solvable GBS groups. In this paper we will always assume that $G$ is a non-solvable GBS group.

A \emph{$G$-tree} is a simplicial tree with an action of $G$ by graph isomorphisms. It can be endowed with a metric, in which case we demand that the action be by isometries.

The set of vertices of a tree $T$ will be denoted by $V(T)$ and the set of edges by $E(T)$. The initial vertex of an edge $e$ is $o(e)$ and its terminal edge is $t(e)$. The opposite edge of $e$ is $\bar e$.

We identify every tree $T$ with its geometric realization. 
A \emph{path} in a tree $T$ is a continuous map from an interval to $T$. We will assume paths are linear on edges. We will frequently identify paths with their image in $T$. A \emph{tight} path is a path which has no backtracking. An \emph{edge path} is a path whose image can be described by whole edges $e_1, \dots, e_k$ such that for $1 \leq i < k$, $t(e_i) = o(e_{i+1})$.

Let $T$ be a $G$-tree. A subgroup $H < G$ is \emph{elliptic} in $T$ if it fixes a point. It is \emph{bi-elliptic} if it fixes two distinct points. When $H$ is not elliptic and contains a loxodromic element, there exists a minimal $H$-invariant subtree which we denote by $T_H$.

A \emph{deformation space} introduced by Forester in \cite{ForesterSplittings} is the space of all minimal $G$-trees which have the same elliptic subgroups, up to $G$-equivariant isometry and homothety. 

The \emph{cyclic deformation space} $\D$ is the deformation space of all $G$-trees where elliptic subgroups are all isomorphic to $\Z$. This is well defined. Indeed there exists a $G$-tree $T$ such that all edge and vertex stabilizers in $T$ are infinite cyclic. Define $\D_T$ as the deformation space containing $T$. For non-elementary GBS groups the set of elliptic subgroups (i.e. subgroups which fix at least a point) is independent of the tree and has an algebraic characterization (see \cite{LevittAutomorphisms} for example) so we set $\D := \D_T$. 

Let $\A$ be a family of elliptic subgroups in $\D$, invariant by conjugacy and by passing to a subgroup. The \emph{restricted deformation space} $\D^\A$ is the subspace of trees in $\D$ where all bi-elliptic groups belong to $\A$.

\bigskip

Let $T$ be a $G$-tree. Let $Y$ be a $G$-invariant proper subforest of $T$. Define the equivalence relation $\sim$ on $T$ as the minimal $G$-equivariant equivalence relation such that all points of a connected component of $Y$ are equivalent. The quotient map $\pi_Y : T \to T/\sim$ is the \emph{collapse} of the subforest $Y$. The tree $T/\sim$ belongs to the same deformation space as $T$ if and only if $Y$ does not contain the axis of any loxodromic element of $G$.

A subgraph of the graph of groups is \emph{collapsible} if its pre-image in the universal cover does not contain the axis of a loxodromic element. Equivalently, when collapsing the corresponding subforest, the stabilizers of the new vertices are elliptic subgroups of $G$.

\begin{rema}
 Collapsibility can be checked using labels.  
 To determine if a subgraph $\Gamma_0 \subset \Gamma$ is collapsible, first of all obtain a new graph from $\Gamma_0$ by deleting all edges with a valence 1 vertex carrying the label $\pm 1$. This yields a graph of groups $\Gamma_1$. Then repeat with $\Gamma_1$. 
 Since the number of edges decreases the procedure will eventually stop, yielding either a graph with no edge, or a graph where all valence 1 vertices have labels different from $\pm 1$. The subgraph $\Gamma_0$ is collapsible if and only if the first case happens.
\end{rema}

Let $\D^\A$ be a restricted deformation space, possibly $\D$ if all edge groups are allowed. A \emph{special factor} with respect to $\D^\A$ is a subgroup $H$ of $G$ with the following property (see \cite{PapinWhitehead} for more details).
There exists a collapse map $\pi: S \to \bar S$ such that $S \in \D^\A$ and $H$ is a vertex stabilizer in $\bar S$.

\medskip

A \emph{marked graph of groups} is a pair $(\Gamma, \sigma)$ where $\Gamma$ is a graph of groups and $\sigma: G \to \pi_1(\Gamma)$ is an isomorphism. Bass-Serre theory gives a correspondance between marked graphs of groups and $G$-trees (\cite{SerreArbresAmalgames}).

Marked graphs of groups offer a nice point of view regarding special factors. A subgroup $H < G$ is a special factor if and only if there exists a marked graph of groups $(\Gamma, \sigma)$ and an identification $\sigma : G \to \pi_1(\Gamma)$ such that $H$ identifies to the fundamental group of a subgraph of $\Gamma$.
In this definition, $\Gamma$ is the quotient $S/G$ with $S \in \D^\A$ such that $H$ is a vertex stabilizer in a collapse $\bar S$.

In the rest of the paper we fix $\A$, so we omit to specify it when referring to special factors.

\begin{defi}
 The outer automorphism $\Phi \in \Out(G)$ is \emph{reducible} with respect to $\D^\A$ if there exists a special factor $A \subset G$ with respect to $\D^\A$ such that the conjugacy class of $A$ is invariant by $\Phi$.
 
 The outer automorphism $\Phi$ is called \emph{fully irreducible} if for all $n \in \N$, $\Phi^n$ is not reducible.
\end{defi}

\begin{rema}
 Since the notion of special factors depends on $\A$, so does the notion of reducible automorphisms. Actually if $\A \subset \A'$ then special factors with respect to $\A$ are special factors with respect to $\A'$, so an automorphism which is reducible for $\A$ is also reducible for $\A'$. The strongest notion of fully irreducible automorphism is achieved with $\D^\A = \D$. We have not explored the possible differences caused by the choice of $\A$ yet.
\end{rema}

\subsection{Computations in GBS groups and trees} \label{sec:algos-de-base} 

In the rest of the paper we will need to perform some algorithmic operations on GBS groups and their automorphisms. In this subsection we justify why and how we can do them. The main operations are summed up by Proposition \ref{prop:algos-de-base}.

In order to deal with the group $G$ algorithmically, we consider elements of $G$ as words in a graph of groups. The trees on which $G$ act also admit a description as a set of words. The definitions in this section are standard, though they are given in the special case of GBS groups whose graphs of groups have cyclic edge and vertex groups, hence the simpler relations. A more detailed description is available in \cite{PapinWhitehead}, Section 3.

Let $\Gamma$ be a graph of groups with cyclic edge and vertex groups. Let $(a_v)_{v \in V(\Gamma)}$ be generators for every vertex group $G_v$ of $\Gamma$ and $(a_e)_{e \in E(\Gamma)}$ be a choice of generators for the edge groups. Each inclusion $G_e \hookrightarrow G_{t(e)}$ is given by a nonzero integer $\lambda(\bar e)$ such that $a_e$ is sent to $a_{t(e)}^{\lambda(\bar e)}$. The integers $\lambda(e)$ is the \emph{label} of $e$.

Let $(t_e)_{e \in E(\Gamma)}$. The \emph{Bass group} $B(\Gamma)$ is the group generated by the elements $a_v$ and $t_e$ for $v \in V(\Gamma)$ and $e \in E(\Gamma)$, and with the following relations:
\begin{itemize}
 \item $t_e = t_{\bar e} ^{-1}$ for $e \in E(\Gamma)$
 \item for all $e \in \Gamma$, with labels $\lambda(e)=p$ and $\lambda(\bar e)=q$, we have $a_{o(e)}^{p}= t_e a_{t(e)}^{q} t_e^{-1}$
\end{itemize}

A \emph{path} in the graph of groups is a word $a_{v_0}^{k_0} t_{e_1} a_{v_1}^{k_1} \dots t_{e_n} a_{v_n}^{k_n}$ where $v_{i-1} = o(e_i)$, $v_i = t(e_i)$ for $1 \leq i \leq n$ and $k_i \in \Z$ for $0 \leq i \leq n$. It is a \emph{loop} if $v_0 = v_n$; in that case we say that it is \emph{based} at $v_0$. Paths represent elements of $B(\Gamma)$. The integer $n$ is the \emph{length} of the path.

Fix a base point $v_0 \in V(\Gamma)$.
The \emph{fundamental group} of $\Gamma$ at vertex $v_0$ is the subgroup $\pi_1(\Gamma, v_0) < B(\Gamma)$ whose elements are represented by loops in $\Gamma$ based at $v_0$. The group is independant of the choice of $v_0$ up to isomorphism.

\bigskip

A path $w:=a_{v_0}^{k_0} t_{e_1} a_{v_1}^{k_1} \dots t_{e_n} a_{v_n}^{k_n}$ in $\Gamma$ is \emph{reduced} if for all $1 \leq i \leq n-1$, $e_i = \bar e_{i+1}$ implies that $\lambda(\bar e_i)$ does not divide $k_i$. 

Given a path $w$, one can algorithmically compute a reduced path which represents the same element of $B(\Gamma)$. 
The reduced path is not necessarily unique, however its length is unique.

\bigskip

The \emph{universal cover} $T_{v_0}$ of the graph of groups $\Gamma$ is a tree which can be constructed as follows: the set of vertices is
\[ 
 \tilde V = \{ \text{paths in } \Gamma \text{ with initial vertex } v_0 \}/ \sim
\]

where $\gamma \sim \gamma'$ if $\gamma$ and $\gamma'$ have the same terminal vertex $v_i \in \Gamma$ and $\gamma^{-1} \gamma' \in G_{v_i}$ as an element of $B(\Gamma)$.

The oriented edges of the universal cover are defined as follows: 
\[ 
 \tilde E =  \left \{ (\alpha,a t_e)/ \begin{array}{l}
                                 \alpha \text{ path in } \Gamma \text{ from } v \text{ to } v', e \in \E_{v'}, a \in G_{v'}
                               \end{array}
  \right \}/ \sim
\]
where $\E _{v'}$ is the set of edges with origin $v'$. The equivalence relation $\sim$ is defined by $(\alpha,a t_e) \sim (\alpha',a' t_{e'})$ if and only if 
$e = e'$ and  $ a^{-1} \alpha^{-1} \alpha' a' \in i_{\bar e}(G_e)$. 
The opposite edge of $(\alpha, a t_e)$ is the edge $(\alpha \cdot a t_e, \bar t_e)$.

The group $\pi_1(\Gamma,v)$ acts on $\tilde V$ and $\tilde E$ by left concatenation of the paths. The quotient of $T_{v_0}$ by this action is $\Gamma$.

\begin{prop}\label{prop:algos-de-base}
 Let $G$ be a GBS group given as a finite set of generators and a marking $G \simeq \pi_1(\Gamma, v)$ where $\Gamma$ is a graph of group with cyclic edge and vertex stabilizers and $v$ is a vertex in $\Gamma$.
 
 One can algorithmically solve the following problems:
 \begin{enumerate}[label = (\roman*)]
  \item word problem: given an element $g \in G$ given as a loop in the graph of groups, find out whether $g$ is the identity
  \item given $x \in \tilde V$ (resp. $e \in \tilde E$) and $g \in G$ given as a loop, compute the image $gx$ 
  \item given vertices $x, x' \in \tilde V$ (resp. edges $e, e' \in \tilde E$) given as paths in the graph of groups, decide if $x \sim x'$ (resp. $e \sim e'$), i.e. if the paths $x,x'$ represent the same point of $T_v$
  \item given $x, y$ in $\tilde V$, find out whether $x,y$ belong to the same orbit, and when they do, find $g \in G$ such that $y = gx$
  \item given two pairs $(x,x'), (y,y')$ of elements of $\tilde V$, find out whether there exists $g \in G$ such that $(y,y') = (gx, gx')$
  \item given $x \in \tilde V$ and $n \in \N$, list all edge paths with length $n$ and first vertex $x$
  \item given an edge $e \in E(\Gamma)$, compute the graph of groups $\Gamma'$ obtained from $\Gamma$ by subdividing $e$
 \end{enumerate}
\end{prop}
\begin{proof}
 (i) Using the identification of the set of generators of $G$ with loops in $\Gamma$, one can compute a loop $w$ in $\Gamma$ for $g$. By \cite[5.2, Theorem 11]{SerreArbresAmalgames} a word $w$ represents the trivial word if and only if it can be reduced to the trivial loop $1 =a_{v_0}^0$, or equivalently if every reduced form is the trivial loop. Since computing a reduced loop is algorithmic, the word problem can be solved.
 
 (ii) To compute the image $gx$ of $x \in V(T)$ given as a path in $\Gamma$ based at $v_0$, it suffices to compute a loop for $g$ based at $v_0$ and concatenate it on the left of the path representing $x$.
 
 (iii) Let $\gamma, \gamma'$ be paths in $\Gamma$ representing $x, x'$. Their first vertex is $v_0$ and their last vertex is a vertex $v_i \in V(\Gamma)$. The  concatenation $\gamma ^{-1} \cdot \gamma'$ is a loop based at $v_i$. It represents an element of $G_{v_i} \subset B(\Gamma)$ if and only if there exists a reduced loop for it with length $0$ (or equivalently, if every reduced loop for $\gamma^{-1} \gamma'$ has length $0$).
 
 (iv) The points $x,y$ belong to the same orbit in $T$ if and only if the paths representing them have the same terminal vertex in $V(\Gamma)$. The element $g \in G$ represented by the loop $y x^{-1}$ satisfies $gx = y$.
 
 (v) We can assume that there exists $g \in G$ such that $gx = y$ and find such a $g$. The set of elements of $G$ such that $gx=y$ is $\{gu, u \in \stab(x)\}$. Let $\gamma_x$ be a path representing $x$ and let $v$ be its terminal vertex in $V(\Gamma)$. By definition of $\tilde V$, $\stab(x) = \{a_x^i, i \in \Z\}$ with $a_x :=\gamma_x a_v \gamma_x^{-1}$. Since $a_v$ is elliptic on $T_{v_0}$ and the tree $T_{v_0}$ is locally finite, the orbit $\langle a_v \rangle \cdot x'$ is finite so there exists $i \geq 1$ such that $a_x^i x' = x'$, and one can algorithmically find such an $i$ by (iii). Thus the pairs $(x, x')$ and $(y,y')$ are in the same orbit if and only if there exists $0 \leq j < i$ such that $g a_x^j x'= y'$. This can be checked algorithmically.
 
 (vi) This can be done by induction. For $n=1$ it suffices to list all edges with origin $x$: if $\gamma_x$ is a path in $\Gamma$ representing $x$ and $v_x$ is the image of $x$ in $V(\Gamma)$, then the edges are 
 \[
  \{(\gamma_x, a_{v_x}^i t_e, e \in E(\Gamma) \text{ s.t. } o(e) = v_x, 0 \leq i < \lambda(e) \}
 \]
 For $n>1$, list all edge paths with length $n-1$. The edge paths with length $n$ are the paths $\tilde \gamma \cdot \tilde e$ where $\tilde \gamma$ is an edge path with length $n-1$ starting with $x$, $\tilde e \in \tilde E$ is an edge whose origin is the last vertex of $\tilde \gamma$, and the last edge of $\tilde \gamma$ and $\tilde e$ are not opposite.
 
 (vii) The subdivision consists in replacing $e$ by edges $e', e''$ with $o(e') = o(e)$, $t(e'')= t(e)$, and $t(e') = o(e'')$ is a new vertex. Labels are the following: $\lambda(e')= \lambda(e)$, $\lambda(\bar e') = \lambda(e'') = 1$, $\lambda(\bar e'') = \lambda(\bar e)$. The new marking must also be computed; the operation consists in replacing $t_e$ in paths in $\Gamma$ by $t_{e'} 1 t_{e''}$.
\end{proof}

A group automorphism $\phi : G \to G$ can be described using a finite set of generators for $G$. We also assume that we know a description of the inverse of the automorphism $G \to \pi_1(\Gamma)$ in terms of the generators of $\pi_1(\Gamma)$ described above. Thus one can compute the image of any element of $\pi_1(\Gamma)$ by $\phi$ and the image is a loop in $\Gamma$.

Maps between trees can be described using the same tool. For a $G$-equivariant map between trees (or a $\phi$-equivariant map; the $G$-equivariant case is the case $\phi = \id$) sending vertex to vertex, it suffices to give the image of every orbit of vertex and edge in the tree. Let $f: T \to T$ be a $G$-equivariant map and let $\mathcal V \subset V(T)$ be a set of representatives for $V(T/G)$, let $\mathcal E \subset E(T)$ be a set of representatives for $E(T/G)$:
\begin{itemize}
 \item the image of a vertex $v \in \mathcal V$ is a path $f(v) \in \tilde V$
 \item the image of an edge $e \in \mathcal E$ is a path in the universal cover: it may be a single vertex or several edges.
\end{itemize}
For another vertex $gv$ (resp. edge $ge$) with $g \in G$ described as a loop in $\pi_1(\Gamma)$, the image is $\phi(g) f(v)$  (resp. $\phi(g) f(e)$) where $\phi(g)$ and $f(v)$ are paths which can be computed separately.

\subsection{Irreducible automorphisms and train track representatives}

Let $T, T'$ be metric $G$-trees. Let $d$ be the distance on $T$ and $d'$ be the distance on $T'$. Let $f : T \to T'$ be a map sending vertex to vertex and edge to non-backtracking edge path. 
The \emph{Lipschitz constant} of $f$ is $\Lip(f) := \sup_{x \neq y \in T} \frac{d'(f(x),f(y))}{d(x,y)}$.

In the rest of the paper, we will always assume that maps between trees send vertex to vertex and send edges to non-backtracking edge paths. We need not assume that maps are linear on edges but for $e \in E(T)$ and $e' \in f(e)$ we will assume that it is linear on $f_{|\{e\}}^{-1}(\{e'\})$.

\medskip

Let $v \in T$ be a vertex. A \emph{turn} at $v$ is an unordered pair of distinct edges $\{e, e'\}$ with origin $v$. If $e=e'$ we call it a \emph{degenerate turn}. A non-backtracking path $\gamma$ \emph{crosses} a turn $\{e, e'\}$ if $e$ and $e'$ appear in $\gamma$.

Fix a map $f : T \to T$. A turn $\{e, e'\}$ is \emph{illegal} if there exists $n \geq 1$ such that $f^n(e)$ and $f^n(e')$ are paths with a common prefix of nonzero length. Otherwise the turn is \emph{legal}. A non-backtracking path $\gamma$ is a \emph{legal path} if every turn crossed by $\gamma$ is legal.

\begin{defi}
We say that $f$ is a \emph{train track map} if $f$ sends every edge $e \in E(T)$ to a legal path. It is a \emph{metric train track map} if in addition the stretch factor on every edge is uniform and equal to $\Lip(f)$, i.e. $\len(f(e)) = \Lip(f) \len(e)$.
\end{defi}

 When $f$ is a train track map, for every $e \in E(T)$ and every $n \in \N$ the path $f^n(e)$ is a geodesic.

 The \emph{bounded cancellation constant} $\BCC(f)$ is a constant introduced in \cite{CooperAutomorphisms} in the case of free groups, with the following property. Let $\alpha, \beta$ be legal paths and let $\alpha \cdot \beta$ be their concatenation. The path $\alpha \cdot \beta$ may not be legal and subsegments of $f(\alpha)$ and $f(\beta)$ may be equal. However the length of the common subsegments is bounded by $\BCC(f)$. The constant exists for the following reason. A piecewise linear map between $G$-trees in the same deformation space is a quasi-isometry (\cite[Remark 3.9]{GuirardelLevitt07}). Thus there is a constant $C$ such that for all $x,y \in T$, $f(x)=f(y) \Rightarrow d(x,y) < C$,  and the constant $\Lip(f) C$ has the property above.
 
\bigskip

Let $\phi \in \Aut(G)$. We say that a map $f: T \to T$ \emph{represents} $\phi$ if it is
$\phi$-equivariant, i.e. for every $t \in T$ and $g \in G$ we have $f(g \cdot t) = \phi(g) \cdot T$.

We call $f: T \to T$ a \emph{train track representative} of $\phi$ if it represents $\phi$ and is a train track map. A train track representative for $\Phi \in \Out(G)$ is a train track representative for some $\phi \in \Phi$.

\begin{rema}
 For all $\Phi \in \Out(G)$, for all $T \in \D$, there exist a $\phi$-equivariant map $f: T \to T$ with $\phi \in \Phi$ sending vertex to vertex and piecewise linear on edges.
 However the train track condition is restrictive.

 A theorem by Bestvina and Handel in \cite{BestvinaHandelTrainTracks} states that every fully irreducible automorphism of $F_N$ has a train track representative. It is not known whether an analogue is true for all GBS groups. Given an arbitrary representative for $\Phi$ one can apply a procedure to try to obtain a train track map but the termination has not been proved, unless $\D$ is finite dimensional (\cite{MeinertTheLipschitzMetric}).
\end{rema}

The \emph{transition matrix} $A(f)$ associated to $f : T \to T$ is a matrix defined as follows. By assumption $f$ sends vertex to vertex and  edge to edge path. Suppose $E(T/G)$ contains $n$ unoriented edges $e_1, \dots, e_n$. The matrix $A(f)$ is the square matrix of size $n$ such that $A(f)_{i,j}$ is the number of edges in the orbit $e_i$ which appear in $f(e_j)$, without taking orientation into account.

A subforest of a tree $T$ is \emph{proper} if it is not $T$. It is \emph{essential} if it contains the axis of an element of $G$.

\begin{lem} \label{lem:FS-sousforet}
 Let $\phi \in \Aut(G)$. The following assertions are equivalent:
 \begin{itemize}
  \item there exists a proper special factor $H <G$ such that the conjugacy class of $H$ is invariant by $\phi$
  \item there exists a tree $T$ and a map $f : T \to T$ representing $\phi$ such that $T$ contains a proper $G$-invariant essential subforest $Y$ with $f(Y) \subset Y$.
 \end{itemize}
\end{lem}
\begin{proof}
 Suppose there exists a proper special factor $H < G$ such that  $H$ is invariant 
 by $\phi$. By definition of special factors there exists a collapse $\pi: T \to \bar T$ with $T \in \D$ such that $H$ is the stabilizer of a vertex $v \in \bar T$. Then $Y:=\pi^{-1}(G \cdot v)$ is a proper $G$-invariant subforest of $T$. Moreover it contains the axis of every loxodromic element of $H$ so it is essential.
 
 Now let us construct the map $f : T \to T$. The subgroup $H$ is a GBS group and $\phi_{|H}$ is an automorphism of $H$. There exists a map $f_H : T_H \to T_H$ representing $\phi_{|H}$. Here is a construction of the map $f_H$. Let $u$ be a vertex in $T_H$. In the collapse $T \to \bar T$, the vertex $u$ is sent to $v$ with $G_v = H$ so $G_u \subset H$. Since $H$ is $\phi$-invariant we also have $\phi(G_u) \subset H$. The point is that there exists a vertex $w \in T_H$ such that $\phi(G_u)$ fixes $w$: there exists $w \in T$ such that $\phi(G_u) \subset G_w$, and by replacing $w$ by its projection on $T_H$ we still have the inclusion.
 
 Define $f_H(gu) = \phi(g) w$ for every $g \in G$. Since there are finitely many of $G$-orbits of vertices in $T_H$, one can repeat this procedure until $f_H$ is fully defined on vertices of $T_H$ and extend to edges by linearity. 
 
 Now we want to extend $f_H$ to $T$. It suffices to define the image of one vertex in every orbit of vertices. Let $v$ be a vertex in $T \setminus G \cdot H$. There exists a vertex $w \in T$ such that $\phi_{|H}(G_v) \subset G_w$. Define $f(v)=w$ and extend by equivariance on $G \cdot v$ by
 \[
  f(gv)= \phi(g) f(v)
 \]
By repeating this for every orbit of vertices and extending linearly on edges, we define $f$ such that $f_{|T_H} = f_H$. Therefore $f(Y) \subset Y$.
 
 \medskip
 
 Conversely suppose there exists $T \in \D$ containing a proper $G$-invariant essential subforest $Y$. Let $Y_0$ be a connected component of $Y$ containing the axis of a loxodromic element $h$. By collapsing the forest $Y$ we obtain a tree $\bar T$ in which $h$ fixes a point $v$. The stabilizer of $v$ is a proper special factor $H$. There are finitely many orbits of vertices in $\bar T$ so there are finitely many conjugacy class of non-cyclic vertex stabilizers in $\bar T$. Since $f(Y) \subset Y$, $f$ induces a map $\bar f$ on $\bar T$. The vertex $v$ must be sent to a vertex with non-cyclic stabilizer. Thus there exists $n \geq 1, k \geq 1$ and $g \in G$ such that $f^{n+k}(v)= gf^k(v)$ so $\phi^{n+k}(H) \subset g \phi^k(H) g^{-1}$, i.e.
 \[
  \phi^{n}(H) \subset g' H g'^{-1}
 \]
 with $g' = \phi^{-k}g$.
By repeating this we obtain a sequence of decreasing special factors \[H \supset g_1^{-1} \phi^n(H) g_1 \supset \dots \supset g_N^{-1} \phi^{nN}(H) g_N \supset \dots\]
 By Corollary 4.11 of \cite{PapinWhitehead} the sequence must be stationary so $H$ and $\phi^n(H)$ are actually conjugate.
\end{proof}

A square non-negative matrix $A$ of size $m$ is \emph{irreducible} if for every $1 \leq i, j \leq m$ there exists $n \in \N$ such that $(A^n)_{i,j} > 0$. The matrix $A$ is \emph{primitive} if there exists $n > 0$ such that all coefficients of $A^n$ are positive. 

A well-known result about primitive and irreducible matrices is:
\begin{theo}[Perron-Frobenius] \label{theo:Perron-Frobenius}
 Let $A$ be a non-negative primitive matrix with size $n\times n$.
 \begin{itemize}
  \item There exists a real eigenvalue $\lambda>0$ (the \emph{Perron-Frobenius eigenvalue}) such that for every other eigenvalue $\mu \neq \lambda$ we have $|\mu|< \lambda$.
  \item The eigenvectors for $\lambda$ are unique up to scalar multiplication and there exists an eigenvector $v$ for $\lambda$ such that $v > 0$.
 \end{itemize} 
 If $A$ is irreducible the same hold except that $|\mu|\leq \lambda$.
\end{theo}
A proof of the theorem can be found in \cite[Theorem 1.1]{Seneta}.

\begin{lem}\label{lem:existeTTirred} 
 Suppose $\phi$ is a fully irreducible automorphism with a train track map $f: T \to T$.
 When $\phi$ is fully irreducible then there exists a collapse $\pi: T \to T'$ with $T' \in \D$ and a train track map $f': T' \to T'$ with $\pi \circ f = f' \circ \pi$ such that $A(f')$ is an irreducible matrix.
 
 Moreover, if the Perron-Frobenius eigenvalue of $A(f')$ is greater than $1$, then $A(f')$ is primitive, so there exists a power $n$ such that $A(f')^n > 0$.
\end{lem}
In that case $f'$ is an \emph{irreducible} train track map. If $A(f')$ is primitive then we call $f$ \emph{primitive}.

\begin{proof}
 Suppose $A(f)$ is not irreducible. There exists a partition $I \cup J = \{1, \dots, n\}$ such that for every $i \in I$, $j \in J$, $A(f)_{i,j}=0$. Then let $\Gamma_J$ be the subgraph of $\Gamma=T/G$ spanned by edges $e_j$ such that $j \in J$. Then the corresponding subforest $T_J \subset T$ is invariant by $f$. By Lemma \ref{lem:FS-sousforet}, $T_J$ is non essential since $\phi$ has no invariant conjugacy class of special factors. 
 
 Let $\pi : T \to  \bar T = T/\sim_{T_J}$ be the collapse of the non essential invariant forest $T_J$. Since $T_J$ is non essential $\bar T \in \D$.  
 We have $f(T_J) \subset J$ so the map $f$ induces a map $\bar f$ on $\bar T$ such that if $\pi$ is the collapse map $T \to \bar T$ then $\pi \circ f = \bar f \circ \pi$.
 
 If $A(\bar f)$ is irreducible then we are done. Otherwise this process can be iterated by collapsing a non essential $f$-invariant forest of $\bar T$. At each step the number of orbits of edges in $\bar T$ decreases so this stops eventually and we obtain the map $f' : T' \to T'$ such that $A(f')$ is irreducible. 
 
 Let us check that the map $f' : T' \to T'$ is a train track map. Let $\pi$ be the collapse map $T \to T'$. For every $k \in \N$ we have $\pi \circ f^k = f'^k \circ \pi$. If $\gamma$ is a non-backtracking path in $T$ then $\pi(\gamma)$ is also non-backtracking. For every edge $e' \in E(T')$, there is a unique edge $e \in \pi^{-1}(e)$ and $f'^k(e') = \pi \circ f^k(e)$. Since $f$ is train track $f^k(e)$ is geodesic and so is $f'^k(e')$.
 
 \bigskip
 
 Let us prove that $A(f')$ is a primitive matrix. Since $f'$ is a train track map, for every $n \in \N$, $A(f'^n) = A(f')^n$. As in \cite[Section 1.4]{Seneta} an irreducible matrix $A$ admits a partition of indices such that $A$ induces a permutation of the classes of the permutation. In terms of tree maps, it means that there exists $k \geq 1$ and a $G$-invariant partition $T_1 \sqcup \dots \sqcup T_k$ of $E(T)$ such that $f'$ induces a permutation of the subforests $T_i, i \in \{1, \dots, k\}$. Moreover $k = 1$ if and only if $A(f')$ is primitive.
 
 By contradiction, suppose $k > 1$.
 There is a power $f'^m$ such that for every $i \in \{1, \dots, k\}$ we have $f'^m(T_i) \subset T_i$. Lemma \ref{lem:FS-sousforet} implies that $T_i$ must be a non-essential, thus collapsible, subforest. Its connected components must be uniformly bounded. 
 
 The Perron-Frobenius of $A(f')$ is greater than $1$ so there is an edge $e$ in $T$ such that the diameter of $f^{mn}(e)$ is unbounded when $n$ goes to infinity. There is a subforest $T_i$ such that $e \in T_i$, and $f^{mn}(e)$ must also be contained in a connected component of $T_i$. This is a contradiction and it implies that $k=1$ and $A(f')$ is primitive.
\end{proof}

Observe that the construction of $f': T' \to T'$ starting from $f: T \to T$ is algorithmic.

\begin{rema} \label{rema:perronFrobenius}
Let $f$ be a train track representative for $\phi$ (i.e. sending edges to legal paths but not necessarily at a uniform speed).
Suppose $A(f)$ is irreducible.
Let $\lambda$ be the Perron-Frobenius eigenvalue for $A(f)$ and let $v$ be the right Perron-Frobenius eigenvector normalized such that its coordinates add up to $1$. The metric on $T$ can be redefined such that $\len(e_i)$ is the $i$-th coordinate of $v$. Then the Lipschitz constant of the new map $f : T \to T$ is uniform on edges and is equal to $\lambda$, so $f$ is a metric train track map.
\end{rema}

Because of Remark \ref{rema:perronFrobenius} all train track maps will be considered as metric train track maps. This point of view is not necessary but it gives some intuition on the behaviour of the iterate images of edges by $f$.

\begin{remas} \label{rem:twist} 
 \begin{itemize}
  \item If $G$ is solvable then it has no special factor. Baumslag-Solitar groups have no fully irreducible automorphism: if $q = pn$ for some $n \in \N$, then $\BS(p,q)$ has an $\Aut(G)$-invariant conjugacy class of special factor (\cite{BouetteThese}). For other Baumslag-Solitar groups or an amalgamated product $\Z *_{\Z} \Z$ all automorphisms have finite order in $\Out(G)$ (\cite{LevittAutomorphisms}). Thus we do not consider these groups when studying fully irreducible automorphisms. 
  \item Suppose $G$ is not solvable and is neither a Baumslag-Solitar group $\BS(p,q)$ nor an amalgamated product $\Z *_{\Z} \Z$: then graphs of groups for $G$ have more than one edge and have a proper subgraph representing a proper special factor.
 Let $f : T \to T$ be an irreducible train track map representing $\phi \in \Aut(G)$. If the Perron-Frobenius eigenvalue of $A(f)$ is $1$ then $T = T \cdot \phi$ and $f$ is an isometry. This implies that $f$ preserves a subgraph of $T/G$ so $\phi$ is not fully irreducible. 
 \end{itemize}
\end{remas}

\begin{cor}\label{coro:construireTTirred}
 Suppose $G$ is a non-solvable Baumslag-Solitar group which is neither a Baumslag-Solitar group nor an amalgamated product $\Z *_{\Z} \Z$. Let $\Phi \in \Out(G)$ with a train track representative. 
 
 There is an algorithm which takes a train track representative $f: T \to T$ for $\Phi$ and
 either computes a primitive train track representative for $\Phi$ or gives a proof that $\Phi$ is reducible.
\end{cor}
\begin{proof}
 The construction of an irreducible representative for $\Phi$ follows Lemma \ref{lem:existeTTirred}. Lemma \ref{lem:existeTTirred} also states that if $\Phi$ is fully irreducible and the Perron-Frobenius eigenvalue is strictly greater than $1$ then the transition matrix is also primitive. By Remark \ref{rem:twist} 2. if the Perron-Frobenius eigenvalue is $1$ then $\Phi$ is reducible. In that case the matrix is not primitive. Since there exists a power $n$ depending on the size of the matrix $A$ such that if $A$ is primitive then $A^n > 0$, one can algorithmically test the primitivity of the transition matrix.
\end{proof}

\subsection{The stable lamination associated to a train track map}

Let $T$ be a $G$-tree in $\D$. Let $\partial T$ be the space of ends of $T$, i.e. the set of equivalence classes of infinite geodesic rays of $T$, where rays $\rho, \rho'$ are equivalent if the Hausdorff distance $d(\rho,\rho')$ is finite. We endow it with the standard topology: a basis of neighbourhoods for $\xi \in \partial T$ is $\{\mathcal V_x, x \in V(T)\}$ where $\mathcal V_x$ is the connected component of $T \setminus \{x\}$ containing a ray for $\xi$.

For any other tree $T' \in \D$ there exists a $G$-equivariant quasi-isometry $f : T \to T'$ which induces a $G$-equivariant homeomorphism $\partial f : \partial T \to \partial T'$. The homeomorphism $\partial f$ does not depend on the choice of $f$ so there is a canonical identification of $\partial T'$ with $\partial T$ for any $T' \in \D$.

A \emph{lamination} of $G$ is a closed, symmetric, $G$-invariant subset of $\partial T \times \partial T \setminus \Delta$, for any tree $T \in \D$, where $\Delta$ is the diagonal. As discussed above, its definition does not depend on the tree $T$.

There is a canonical action of $\Aut(G)$ on $\partial T$ and thus on the set of laminations: let $T \in \D$ and $\phi \in \Aut(G)$ and let $f: T \to T$ be a representative for $\phi$. It induces a homeomorphism of $\partial T$ which does not depend on the choice of $f$.

Since laminations are $G$-invariant, the action of $\Aut(G)$ on $\partial T$ yields an action of $\Out(G)$ on the set of laminations.

Let $\Lambda$ be a lamination. For $T \in \D$, the \emph{realization} $\Lambda_T$ of $\Lambda$ in $T$ is the $G$-invariant set of unoriented bi-infinite geodesics whose endpoints belong to  $\Lambda \subset \partial T \times \partial T \setminus \Delta$. The geodesics of the realization are called \emph{leaves} of $\Lambda$. If $T' \in \D$ and $f : T \to T'$ is a quasi-isometry, then for any leaf $\lambda \in \Lambda_T$, the geodesic obtained by tightening $f(\lambda)$ is a leaf of $\Lambda_{T'}$. Conversely all leaves of $\Lambda_{T'}$ are obtained that way.

A \emph{leaf segment} of $\Lambda_T$ is a segment of $T$ contained in a leaf in $T$.

We now introduce the \emph{stable lamination} of an automorphism with a irreducible train track.

\begin{defi}
 Let $\Phi \in \Out(G)$ and $f : T \to T$ be an irreducible train track representative for $\Phi$.
 The \emph{stable lamination} $\Lambda^+$ is defined by its realization in $T$. A bi-infinite geodesic $\lambda$ lies in $\Lambda^+_T$ if and only if for every leaf segment $\sigma \subset \lambda$, there exists an edge $e \in E(T)$ and $n \in \N$ such that $\sigma \subset f^n(e)$.
\end{defi}

\begin{rema}
 Since all edges are legal, all leaf segments are also legal, so there is no cancellation in leaves in $T$ when applying $f$. Thus leaves of $\Lambda^+_T$ are legal. Viewing $f$ as a metric train track map, they are also uniformly stretched by the factor $\Lip(f)$. The set $\Lambda^+_T$ is stable by $f$. It implies that $\Lambda^+$ is stable by $\Phi$.
\end{rema}

\begin{rema}
 If $\Phi \in \Out(G)$ admits an irreducible train track representative $f:T \to T$ and thus a stable lamination $\Lambda^+$ can be defined, then $\Phi^n$ admits $f^n :T \to T$ as a train track representative. Then the stable lamination associated to $f^m$ is equal to $\Lambda^+$. Indeed, a subsegment of $f^k(e)$ is also a subsegment of $f^{nm}(e)$ for some $m \in \N$ : since for $N$ big enough, $f^N(e)$ contains a translate of $e$, then for $m$ big enough $f^{nm-k}(e)$ contains a translate of $e$ and $f^{nm}(e)$ contains a translate of $f^k(e)$. 
\end{rema}

\begin{defi} 
 A lamination $\Lambda$ is \emph{minimal} with respect to a tree $T$
 if it satisfies the following condition:   $\forall \lambda, \lambda' \in \Lambda_T$, $\forall I$ leaf segment in $\lambda$, $\exists g \in G$ such that $gI \subset \lambda'$. In other words, all leaves of a minimal lamination have the same leaf segments up to the action of $G$.
\end{defi} 
\begin{rema} \label{rem:minimal}
 Let $T, S \in \D$.
 A lamination $\Lambda$ is minimal with respect to $T$ if and only if it is minimal with respect to $S$. Indeed let $f : T \to S$ be a $G$-equivariant quasi-isometry. There exists $C > 0$ depending on $f$ such that the image of any geodesic path in $T$ is in the $C$-neighbourhood of a geodesic path in $S$. 
 
 Let $\lambda_S, \lambda'_S$ be leaves of $\Lambda_S$. Let $\lambda_T, \lambda_T'$ be the corresponding leaves in $\Lambda_T$. Let $I_S \subset \lambda_S$ be a leaf segment. There exists $I_T \subset \lambda_T$ such that $f(I_T)$ contains a $C$-neighbourhood of $I_S$ in $\lambda_S$. There also exists $g \in G$ such that $gI_T \subset \lambda_T'$. The segment $f(gI_T) = gf(I_T)$ is in a $C$-neighbourhood of $\lambda_S'$ so its central part containing $gI_S$ is contained in $\lambda'_S$.
\end{rema}

\begin{lem}
 The stable lamination associated to an automorphism with an irreducible train track representative is minimal.
\end{lem}
A proof can be found in \cite[Lemma 1.2]{BestvinaFeighnHandelLaminations} for free groups, and it can be adapted to the case of GBS groups.

\begin{defi}
 A leaf $\ell$ of the realization of a lamination in a tree is \emph{quasi-periodic} if for every $l > 0$ there exists $C_l>0$ such that if $\sigma$ is a leaf segment of $\ell$ of length at most $l$ then any leaf segment of $\ell$ of length at least $C_l$ contains a translate of $\sigma$.
\end{defi}
\begin{lem}
 Let $\Phi \in \Out(G)$ with a train track $f : T \to T$. Let $\ell$ be a leaf of $\Lambda^+_T$. The leaf $\ell$ is quasi-periodic.
\end{lem}
\begin{proof}
 Let $l > 0$.
 Let $e \in E(T)$. Since the stable lamination is minimal, there exists $k \in \N$ depending on $l$ such that $f^k(e)$ contains an edge path in every orbit of edge paths with length at most $l$ crossed by leaves of the lamination.
 
 There exists $n \in \N$ such that for every $e' \in E(T)$ and $m \geq n$, the path $f^m(e')$ contains an edge in the orbit of $e$.
 
 Let $C_l:= \max_{e \in E(T)} \len( f^{n+k}(e'))$. Let $\sigma \subset \ell$ be a leaf segment with length at least $C_l$. There exists $m \in \N$ and $e' \in E(T)$ such that $\sigma \subset f^m(e')$. Moreover $m \geq n+k$. The path $f^{m-k}(e')$ contains an edge in the orbit of $e$ so $f^m(e')$ contains a subpath in the orbit of $f^k(e)$, which itself contains a path in every orbit of edge paths with length at most $l$.
\end{proof}

Let us fix an automorphism $\Phi \in \Out(G)$, and suppose it admits a train track representative and thus a stable lamination $\Lambda^+$. Let  $S \in \D$, not necessarily the train track representative. Let $v$ be a vertex in $S$. Let $\lambda$ be a leaf in the stable lamination $\Lambda^+_S$. The \emph{Whitehead graph} of the leaf $\lambda$ at the vertex $v$ is the graph $\Wh_{S}(\lambda, v)$ such that
\begin{itemize}
 \item vertices are edges of $S$ with origin $v$
 \item there is an edge $e - e'$ if there exists $g \in G$ such that $\{e, e'\}$ is a turn in $g\lambda$
\end{itemize}

The vertex stabilizer $G_v$ acts naturally on $\Wh_{S}(\lambda, v)$. The stable lamination is minimal so all leaves of the stable lamination have the same subsegments up to translation. In particular they have the same turns up to translation so for any leaf $\lambda'$ in $\Lambda^+$ we have $\Wh_S(\lambda',v) = \Wh_S(\lambda,v)$. We may as well define the Whitehead graph of the stable lamination at vertex $v$ by $\Wh_S(\Lambda^+,v)= \Wh_S(\lambda,v)$ for any leaf $\lambda$.

\begin{lem} \label{lem:algoFeuilles}
 Let $\Phi \in \Out(G)$ be an automorphism with an irreducible train track representative $f : T \to T$ and associated stable lamination $\Lambda^+$. The Whitehead graphs $\Wh_T(\Lambda^+, v)$ can be computed algorithmically.
\end{lem}
\begin{proof}
 In order to compute Whitehead graphs we need to find all orbits of turns taken by the lamination. 
 
 Let $e \in E(T)$. Since $f$ is irreducible, every leaf segment of $\Lambda^+_T$ is a subsegment of a translate of $f^n(e)$ for some $n \in \N$. Therefore all turns taken by $\Lambda^+_T$ appear in $f^n(e)$ for some $n \in \N$.
 
 There are two ways a turn can arise in some $f^n(e)$. If $\{e_1, e_2\}$ is a turn in $f^n(e)$ then either it is a turn in $f(e')$ for some $e' \in E(T)$ (first type), or it is the turn between $f^k(e_1')$ and $f^k(e_2')$ where $\{e_1', e_2'\}$ is a turn of the first type (second type).
 
 There exists a power $f^n$ depending only on the size of $A(f)$ such that for every edge $e$ of $T$, $f^n(e)$ crosses every orbit of edge in $T$. Therefore all turns of the first type appear in $f^{n+1}(e)$. Let $K$ be the number of orbits of turns in $T$, which is finite. Then all turns of the second type appear in $f^{n+1+k}(e)$ for some $k \leq K$. This gives a bound on the number of iterations needed to find all turns.
\end{proof}

\bigskip

The following results may apply to more general laminations.
Let $\Lambda$ be a lamination. We say that a finitely generated subgroup $H \subset G$ containing a hyperbolic element \emph{carries} $\Lambda$ if there exists $S \in \D$ such that any leaf in $\Lambda_S$ is contained in a translate of the minimal subtree for $H$ in $S$. The condition only depends on the conjugacy class of $H$.
Moreover it does not depend on a choice of $S$:
\begin{lem} \label{lem:changerArbreFeuille}
 If a finitely generated subgroup $H$ carries $\Lambda$ then in any $S \in \D$, any leaf of $\Lambda_S$ is contained in a translate of the minimal subtree $S_H$.
\end{lem}
\begin{proof}
 Let $S \in \D$. A bi-infinite geodesic $\lambda$ is contained in $S_H$ if and only if it lies in a bounded neighbourhood of $S_H$.
 
 Since $H$ carries $\Lambda$, there exists a tree $S' \in \D$ such that all leaves of $\Lambda_{S'}$ are contained in $S'_H$.
 
 Let $f: S' \to S$ be a $G$-equivariant map. It is a quasi-isometry. Let $\lambda$ be a leaf of $\Lambda_{S'}$. Up to translating it we may suppose that it is contained in $S'_H$. Then $f(\lambda)$ is contained in $f(S'_H)$. By minimality the tree $S_H$ is contained in $f(S'_H)$. The diameter of $S'_H/H$ is bounded, also by minimality, so the diameter of $f(S'_H)/H$ is also bounded. Therefore $f(S'_H)$ must be contained in a $c$-neighbourhood of $S_H$ for some $c>0$.
 
 The image $f(\lambda)$ must then be contained in a $c$-neighbourhood of $S_H$. The tightened leaf $[f(\lambda)]$ is a geodesic contained in $f(\lambda)$ so it is contained in $S_H$.
\end{proof}

For minimal laminations such as the stable lamination of an automorphism, we have the following:
\begin{prop} \label{prop:feuilleLamination}
 Let $\Lambda$ be a non-empty minimal lamination. Let $H$ be a special factor of $G$. Then $\Lambda$ is carried by $H$ if and only if there exists $T \in \D$ such that there exists a leaf of $\Lambda_T$ contained in the minimal subtree for $H$ in $T$.
\end{prop}
 
 \begin{proof}
 The direct implication is immediate. Conversely, by Lemma \ref{lem:changerArbreFeuille}, if there exists $T \in \D$ and a leaf $\lambda \in \Lambda_T$ contained in the minimal subtree $T_H$, then for all $T' \in \D$, the realization of $\lambda$ in $\Lambda_{T'}$ is also contained in $T'_H$.
 
 We choose $T' \in \D$ such that $T'_H$ and its translates are disjoint: for all $g \in G$, $g T'_H \cap T_H \neq \varnothing \Rightarrow g \in H$. The leaf $\lambda$ is contained in $T'_H$. If $\lambda'$ is another leaf, its segments are translates of segments of $\lambda$. This implies that $\lambda'$ all edges crossed by $\lambda'$ belong to $G \cdot T'_H$. As translates of $T'_H$ are disjoint $\lambda'$ must be contained in a single translate of $T'_H$.
 \end{proof}

\begin{lem}\label{lem:malnormal}
 Let $A$ be a special factor for $G$. Let $T \in \D$. Let $T_A$  be the minimal subtree for $A$. There exists $C>0$ such that for all $g \in G \setminus A$ we have $\diam(T_A \cap g \cdot T_A) \leq C$.
\end{lem}
\begin{proof}
By definition of a special factor, there exists a tree $S \in \D^\A$ such that the translates of the minimal subtree $S_A$ are disjoint.
 
 There exists a $G$-equivariant application $f : S \to T$ which is a quasi-isometry and such that $T_A \subset f(S_A)$. There exists a constant $C$ depending on $f$ such that if two subtrees are disjoint, the diametre of their intersection is bounded by $C$.
\end{proof}
 
We deduce the following, which can be applied to the stable lamination of an automorphism in our context:
\begin{prop} \label{prop:factMinLamination}
 Let $\Lambda$ be a minimal lamination. There exists a unique special factor relative to $\D$ and minimal for inclusion which carries $\Lambda$.
\end{prop}
\begin{proof}
 The existence is a consequence of the descending chain condition given in Corollary 4.11 of \cite{PapinWhitehead}. For uniqueness we need to show that if a leaf of $\Lambda$ is contained in the minimal subtree for two special factors $A$ and $A'$ then it is also contained in the minimal subtree for $A \cap A'$.
 
 It suffices to prove that the intersection $T_A \cap T_{A'}$ is contained in a bounded neighbourhood of $T_{A \cap A'}$. By Lemma \ref{lem:malnormal} there exists $C > 0$ such that for every $g \in G$, $\diam(T_A \cap g T_A) \geq C \Rightarrow g \in A$ and $\diam( T_{A'} \cap g T_{A'} ) \geq C \Rightarrow  g \in A'$.
 
 There exists $J > 0$ depending on $C$ and $T$ such that every segment $\sigma$ of $T$ with length at least $J$ contains at least two segments of length $C$ in the same $G$-orbit and with same orientation in $\sigma$.
 
 Let $\sigma \subset T_A \cap T_{A'}$ be a segment with length at least $J$. There exists a segment $I \subset \sigma$ with length $C$ and $g \in G$ such that $gI \subset \sigma$. Then the diameter of $T_A \cap gT_A$ is greater than $C$ so $g \in A$. Similarly $g \in A'$. The axis of $g$ is contained in $T_{A \cap A'}$. Thus $T_A \cap T_{A'}$ is contained in a $J$-neighbourhood of $T_{A \cap A'}$.
\end{proof}

 \section{Whitehead graphs of the lamination and reducibility} \label{sec:detecter}

 In order to decide whether an automorphism with a train track representative is fully irreducible, the Whitehead graphs of the stable lamination give important information. In order to state the main theorem we need to introduce first pseudo-periodic conjugacy classes, which will be developed further in Section \ref{sec:nielsen}.
 \begin{defi}
  Let $\Phi \in \Out(G)$ and let $\phi \in \Phi$.
 
  The conjugacy class of an element $g \in G$ is \emph{pseudo-periodic} for $\Phi$ if $\| \phi^n(g)\|_T$ is bounded.
 \end{defi}
 A fact worth mentioning is that the minimal special factor containing a pseudo-periodic element is periodic. We prove this in Section \ref{sec:nielsen}.
 
 The aim of this section is to prove:
\begin{theo} \label{theo:critere2}
 Let $\Phi$ be an automorphism of $G$ with a train track $f: T \to T$ with no simple pseudo-periodic conjugacy class. Then $\Phi$ is fully irreducible if and only if for every $v \in T$, the  Whitehead graph $\Wh_T(\Lambda^+_f,v)$ is connected.
\end{theo}

 Propositions \ref{prop:WhNonConnexeImpliqueRed} and \ref{prop:WhCnxIrr} are GBS equivalents for analogue results by Kapovich (\cite[Proposition 4.1, Proposition 4.2]{KapovichAlgorithmic}). Along with Corollary \ref{coro:feuille} they prove Theorem \ref{theo:critere2}.
 
\begin{prop}\label{prop:WhNonConnexeImpliqueRed}
 Let $\Phi \in \Aut(G)$ be an automorphism with an irreducible train track representative $f : T \to T$ and an associated stable lamination $\Lambda^+$.
 
 Let $S \in \D$ be any tree. Let $f' : S \to S$ be a representative for $\Phi$. If there exists $v \in S$ such that the Whitehead graph $\Wh_S(v, \Lambda^+)$ is not connected and such that the stabilizer of some of its connected components is in $\A$, then $\Phi$ is reducible in $\D^\A$.
\end{prop}

\begin{rema}
 If $\D^\A = \D$ the condition on stabilizers of connected components is always true.
\end{rema}

\begin{lem} \label{lem:feuillePortee}
 Let $\Phi$ be an outer automorphism having an irreducible train track representative. Let $S$ be a $G$-tree. If there exists a special factor whose minimal subtree in $S$ contains a leaf of the stable lamination, then $\Phi$ is reducible.
\end{lem}
\begin{proof}
 Let $\lambda$ be a leaf of the stable lamination. There exists a unique minimal special factor $H$ carrying the leaf $\lambda$ (Proposition \ref{prop:factMinLamination}). By Proposition \ref{prop:feuilleLamination}, every leaf of the lamination is contained by a translate of $H$. Let $\phi \in \Phi$. Since the stable lamination is $\phi$-invariant, all leaves are also carried by translates of $\phi(H)$. By minimality there is an element $h$ such that we have $H \subset h\phi(H)h^{-1}$. It follows that $H=h \phi(H) h^{-1}$. Otherwise, construct a decreasing sequence
 \[
   \dots \subset h_n^{-1} \phi^{-n}(H) h^n \subset \dots \subset H
 \]
and by \cite[Corollary 4.11]{PapinWhitehead} this sequence is stationary. This implies that $H$ is conjugate to $\phi(H)$. 
 \end{proof}

\begin{proof}[Proof of Proposition \ref{prop:WhNonConnexeImpliqueRed}.]
 Let $\phi \in \Aut(G)$ with a train track representative $f: T \to T$. Let $f': S \to S$ be a representative for $\phi$. Suppose there exists a vertex $v$ such that the Whitehead graph $\Wh_S(v,\Lambda^+)$ is non-connected with at least one connected component with stabilizer in $\A$. We will construct a tree $S'$ such that there exists a collapse $S' \to S$, and such that the lifts of the leaves of $\Lambda^+$ in $S'$ avoid an orbit of edges.

 Let $C_1, \dots, C_k$ be the connected components of $\Wh_S(v,\Lambda^+)$.
 
 Let $E$ be the star with $k$ edges $e_1, \dots, e_k$. Replace $v$ by $E$ by attaching the end of $e_i$ to edges in $C_i$. By extending this construction by $G$-equivariance, 
 we define an expansion at vertex $v$.
 
 The edge $e_i$ has the same stabilizer as the corresponding connected component $C_i$. It is possible that the construction created an edge with stabilizer not in $\A$. Collapse these edges. By assumption, at least one edge in $E$ has stabilizer in $\A$, so not all edges are collapsed.
 
 We obtain a tree $S'$ with a collapse map $\pi : S' \to S$.
 The leaves of the stable lamination can be lifted in $S'$. 
 
 By construction leaves in $S'$ do not cross the edges which are collapsed by $\pi$. Thus there is a special factor $H$ which carries a leaf of the lamination. By Lemma \ref{lem:feuillePortee}, $\phi$ is reducible.
\end{proof}

\begin{prop}\label{prop:WhCnxIrr}
 Let $\Phi \in \Out(G)$. Let $f:T \to T$ be a train track representative for $\Phi$. Suppose that the incidence matrix $A(f)$ is primitive. Let $\Lambda^+$ be the stable lamination.  
 
 Then if for every $v \in T$, the Whitehead graph $\Wh_T(\Lambda^+, v)$ is either connected, or disconnected such that no connected component has a stabilizer in $\A$, then no leaf of the lamination is carried by a proper special factor.
\end{prop}

\begin{proof}[Proof of Proposition \ref{prop:WhCnxIrr}.]
 First let us prove the proposition under simplified assumptions: we suppose that all Whitehead graphs are connected.
 
 By contradiction, let $\Phi, f, T$ be as in the proposition and 
 such that all Whitehead graphs of the stable lamination are connected. Let $\phi \in \Phi$ be such that $f$ is $\phi$-equivariant. Let $A$ be a special factor such that the leaf $\lambda$ is contained in the minimal subtree $T_A$ of $A$ in $T$. By Proposition \ref{prop:feuilleLamination} this implies that $A$ carries the stable lamination and that every leaf of it is contained in a translate of $T_A$.
 
 Let $e$ be an edge in $T_A$. Let $e'$ be any edge of $T$. Up to reversing the orientation of these edges, there exists a geodesic path $e = e_0 \dots e_n= e'$ (see Figure \ref{fig:recouvParFeuilles} for what follows). Turns in this path are not necessarily crossed by any leaf of $\Lambda^+$. However, Whitehead graphs are connected. At every turn $\bar e_i, e_{i+1}$ at vertex $v_i$, there exist edges $\bar e_i = \epsilon_0,\dots, \epsilon_k= e_{i+1}$ with origin $v_i$ such that $\epsilon_j, \epsilon_{j+1}$ is a turn crossed by a leaf. Thus $\bar \epsilon_j \epsilon_{j+1}$ is a leaf segment.
 
 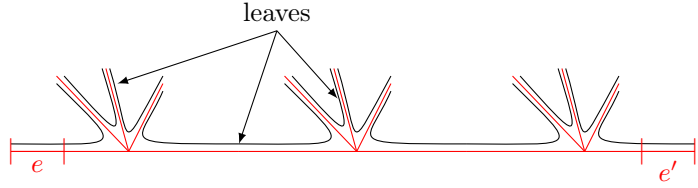
\begin{figure}[h] 
  \centering
   \begin{tikzpicture}
\draw (-0.5,0) .. controls (1,0) and (1,-0.1) .. (0.1,0.8);
\draw (0.7,1) .. controls (1,0) and (1,0) .. (0.2,0.9);
\draw (0.8,1) .. controls (1.1,0) and (0.9,-0.2) .. (1.5,0.9);
\draw (2.5,0) .. controls (0.8,0) and (1.2,0) .. (1.5,0.7);

\draw [color=red] (-0.5,-0.1) --(2.5,-0.1);
\draw [color=red] (1.05,-0.1) --(0.1,0.9);
\draw [color=red] (1.05,-0.1) --(0.75,1);
\draw [color=red](1.05,-0.1) --(1.5,0.8);

\draw [color=red] (-0.5,-0.1)node {$|$} --(0.2,-0.1)node {$|$} node [midway, below] {$e$};

\draw [-latex] (3,1.5)node [above] {leaves} -- (2.5,0);
\draw [-latex] (3,1.5) -- (0.9,0.8);
\draw [-latex] (3,1.5) -- (3.8,0.6);

\begin{scope}[xshift=3cm]
\draw (-0.5,0) .. controls (1,0) and (1,-0.1) .. (0.1,0.8);
\draw (0.7,1) .. controls (1,0) and (1,0) .. (0.2,0.9);
\draw (0.8,1) .. controls (1.1,0) and (0.9,-0.2) .. (1.5,0.9);
\draw (2.5,0) .. controls (0.8,0) and (1.2,0) .. (1.5,0.7);

\draw [color=red] (-0.5,-0.1) --(2.5,-0.1);
\draw [color=red] (1.05,-0.1) --(0.1,0.9);
\draw [color=red] (1.05,-0.1) --(0.75,1);
\draw [color=red](1.05,-0.1) --(1.5,0.8);

\end{scope}

\begin{scope}[xshift=6cm]
\draw (-0.5,0) .. controls (1,0) and (1,-0.1) .. (0.1,0.8);
\draw (0.7,1) .. controls (1,0) and (1,0) .. (0.2,0.9);
\draw (0.8,1) .. controls (1.1,0) and (0.9,-0.2) .. (1.5,0.9);
\draw (2.5,0) .. controls (0.8,0) and (1.2,0) .. (1.5,0.7);

\draw [color=red] (-0.5,-0.1) --(2.5,-0.1);
\draw [color=red] (1.05,-0.1) --(0.1,0.9);
\draw [color=red] (1.05,-0.1) --(0.75,1);
\draw [color=red](1.05,-0.1) --(1.5,0.8);

\draw [color=red] (1.8,-0.1)node {$|$} --(2.5,-0.1)node {$|$} node [midway, below] {$e'$};
\end{scope}
\end{tikzpicture}
  \caption{Cover the path between $e$ and $e'$ with overlapping leaf segments.}\label{fig:recouvParFeuilles}
  
 \end{figure}
 
  This yields a sequence of leaf segments which cover the path between $e$ and $e'$ and such that two consecutive segments overlap over at least the length of an edge. 
  
  Since leaf segments are legal, they are all stretched uniformly by $f$. There is a power $m$ of $f$ such that the image of any edge by $f^m$ is strictly longer than constant $C$ given by Lemma \ref{lem:malnormal}. 
  
  If we apply $f^m$ to the leaf segments above, we obtain longer leaf segments such that two consecutive segments overlap by a length greater than $C$. Since all leaf segments are contained in translates of $T_A$, they must be contained in the same translate $h_{e'} T_A$. In particular $f^m(e)$ and $f^m(e')$ are both contained in $h_{e'} T_A$, by Lemma \ref{lem:malnormal}.
  
  Since $f^m(e)$ is longer than $C$ and $f^m(e) \subset h_{e'} T_A$ for any $e' \in E(T)$, the element $h_{e'} \in G$ actually does not depend on $e'$ and we denote it by $h$.

  Therefore there exists $m \in \N$ and $h \in G$ such that $f^m(e') \subset h T_A$ for all $e' \in E(T)$. This implies that $f^m(T) \subset T_A$, which gives the desired contradiction since $f^m$ is surjective by minimality of the action.
  
  \bigskip
  
  Now suppose that Whitehead graphs may be disconnected, though stabilizers of connected components do not belong to $\A$. Let $A$ be a special factor such that $T_A$ contains a leaf $\lambda$ of the stable lamination.
  
  Let $\mathcal R$ be the smallest $G$-equivariant equivalence relation on $E(T)$ such that two edges in a turn crossed by a leaf of the lamination belong to the same equivalence class. It defines a partition of $E(T)$. If $e \mathcal R e'$ then all edges on the geodesic path between $e$ and $e'$ also belong to the same equivalence class. Each equivalence class spans a connected subtree of $T$. When all Whitehead graphs are connected, there is a single equivalence class.

  Now we prove that there exists $m \in \N$ such that $f^m(S_{i_0})$ is contained in a translate of $T_A$. The argument is the same as in the simple case of the proof: take $e \in T_A \cap S$ on the leaf $\lambda$. There is $m \in \N$ depending uniquely on $T, A, f$ such that for any $e' \in S$ the images $f^m(e)$ and $f^m(e')$ are in the same translate of $T_A$. This works well since the restriction of any Whitehead graph to $S_{i_0}$ is connected. Then we have $f^m(S_{i_0}) \subset g T_A$ for some $g \in G$.
   
  As $f^m(S) \subset gT_A$ and $f^m(S_{i_0})$ is unbounded, Lemma \ref{lem:malnormal} implies $\phi^m(\stab S) \subset \stab f^m(S_{i_0}) \subset gAg^{-1}$, so $\stab(S) \subset \phi^{-m}(gAg^{-1})$.
   
  \bigskip
  
  Let $\{S_i\}_{i \in I}$ be the set of subtrees of $T$ corresponding to equivalence classes of $\mathcal R$. The intersection $S_i \cap S_{j}$ for $i \neq j$ is either empty or a single vertex. This means that $\{S_i\}_{i \in I}$ is a \emph{transverse covering} of $T$ (see \cite[Def. 4.6]{GuirardelLimitGroups}).
  
  There exists $i_0 \in I$ such that $\lambda$ is contained in $S_{i_0}$. Since the transition matrix of the train track is primitive, $\lambda$ crosses all orbits of edges in $T$ so all subtrees $S_i$ are translates of $S_{i_0}$. Construct $\hat T$ (the \emph{skeleton of the transverse covering}, \cite[Def. 4.8]{GuirardelLimitGroups})
  as follows. It is a bipartite graph. Vertices of $\hat T$ are $s_h$ with $ h \in G/\stab S_{i_0}$ which correspond to translates of $S_{i_0}$ and $x_v$ for $v \in T$ such that $v$ is the intersection of two distinct translates of $S_{i_0}$. There is an edge $s_h - x_v$ if $v \in h S$. This yields a $G$-tree with vertex stabilizers of $x_v$ elliptic in $\D$ and vertex stabilizers of $s_h$ conjugate to $\stab S$ \footnote{This tree may also be obtained by replacing vertices of $T$ by stars on the connected components of their Whitehead graph, as in the proof of Proposition \ref{prop:WhNonConnexeImpliqueRed}, then collapsing copies of $S_{i_0}$.}. The edge group of $s_h - x_v$ is the stabilizer of the component of $\Wh_T(\Lambda^+,v)$ which intersects $hS_{i_0}$, so it is not in $\A$.

  Construct a map $\hat T \to R$ where $R$ is a collapse of an element of $\D^\A$ and has a vertex stabilizer $\phi^{-m}(A)$. Map the vertex $s_h$ of $\hat T$ for $h \in G$ to the vertex of $R$ with stabilizer $h\phi^{-m}(A)h^{-1}$ and map every vertex $x_v$ to a vertex of $R$ whose stabilizer contains $G_{v}$. All edge groups in $R$ are in $\A$ while no edge group in $\hat T$ is in $\A$. Thus every edge of $\hat T$ is sent to a point in $R$. By continuity the image of $\hat T$ is a single point in $R$ and is $G$-invariant. By minimality of the action on $R$, the tree $R$ must be a single point and $\phi^{-m}A=G$ so $A = G$, which is a contradiction.
 \end{proof}

\begin{cor} \label{coro:feuille}
If $\Phi \in \Out(G)$ satisfies the hypotheses of Proposition \ref{prop:WhCnxIrr}, and no simple element in $G$ is pseudo-periodic for $\Phi$, then $\Phi$ is fully irreducible.
\end{cor}

\begin{proof} 
 By contradiction, let $\Phi$ be an automorphism satisfying the hypotheses of the corollary. Let $\phi \in \Phi$. 
 
 Let $H$ be a special factor of $G$. Up to replacing $\Phi$ by $\Phi^k$ and to choosing $\phi \in \Phi^k$, we may assume that $\phi(H)=H$.

 Let $T_H$ be the minimal subtree of $H$. Let $h \in H$ be a loxodromic element whose axis is contained in $T_H$, so $h$ is simple. By assumption $\|\phi^n(h)\|$  goes to infinity.

 We will find a contradiction to Proposition \ref{prop:WhCnxIrr} by finding a leaf of $\Lambda^+$ contained in $T_H$. 
 
 Let $e_1, \dots, e_m$ be the edges of a fundamental domain for $h$. Some fundamental domain of the axis of $\phi^n(h)$ can be written as a concatenation of (maybe empty) sets $J_i^n$ where $J_i^n \subset f^n(e_i)$ for $1 \leq i \leq m$.

 The non-degenerate $J_i^n$ are leaf segments of $\Lambda^+$, some of them may be empty. Since $\|\phi^n(h)\|_T$ goes to infinity, at least one of the $J_i^n$ must also go to infinity when $n$ goes to infinity. Besides we have $J_i^n \subset T_{\phi(H)} = T_H$ for all $n \in \N$.

  Suppose that $(J_{i_0}^n)_{n\in \N}$ tends to infinity. The segments $J_{i_0}^n$ are arbitrarily long leaf segments, all contained in $T_H$. Let $\lambda$ be a leaf of $\Lambda^+$. Then every leaf segment of $\lambda$ is contained in a translate of $J_{i_0}^n$ for some $n \in \N$.
  
  By Lemma \ref{lem:malnormal} there exists $C$ such that if the intersection of two translates of $T_H$ has diameter greater than $C$, then the translates are equal.
  
  Since $\lambda$ is quasi-periodic there exists $L > 0$ such that for every leaf segment $\sigma$ of length $3C$, every leaf segment $\gamma \subset \lambda$ with $\len(\gamma) \geq L$ contains a translate of $\sigma$.
  
  There exists $i_0$ and $n \in \N$ such that $\len(J_{i_0}^n) \geq L$. By minimality of $\Lambda^+$, the leaf $\lambda$ contains a translate of $J_{i_0}^n$ so every leaf segment $\sigma \subset \lambda$ with $\len(\sigma) = 3C$ is a translate of a subsegment of $J_{i_0}^n \subset T_H$. Tile $\lambda$ with such leaf segments of length $C$ such that two consecutive tiles overlap on more than $C$. Each segment of the tiling is contained in a translate of $T_H$ and by Lemma \ref{lem:malnormal} all these translates must be equal. Thus $\lambda \subset g T_H$ for some $g \in G$ and Proposition \ref{prop:feuilleLamination} implies that $H$ carries $\Lambda$.  
\end{proof}

\section{Nielsen paths and pseudo-periodic elements} \label{sec:nielsen}

The aim of the present section is to define and describe \emph{pseudo-periodic} conjugacy classes of an outer automorphism $\Phi \in \Out(G)$. They are analogues of periodic conjugacy classes for an automorphism of a free group.

Pseudo-periodicity is a weaker notion than periodicity: the conjugacy class of $g \in G$ is pseudo-periodic if the translation length $\|\phi^n(g)\|_T$ is bounded for some, equivalently any, tree $T \in \D$, and some $\phi \in \Phi$. While it automatically implies that the conjugacy class of $g$ is periodic in the free group case, it is not always true for GBS groups. 

An outer automorphism $\Phi \in \Out(G)$ is \emph{pseudo-atoroidal} if it has no pseudo-periodic conjugacy classes.

The point is that when there exists a \emph{simple} element $g \in G$ whose conjugacy class is pseudo-periodic, then the conjugacy class of the unique proper minimal special factor containing $g$ is $\Phi$-periodic. In Section \ref{sec:detecter} we gave a criterion for reducibility, such that automorphisms satisfying the criterion are irreducible if and only if they do not have any pseudo-periodic conjugacy classes.

For algorithmic purposes we need to find the pseudo-periodic conjugacy classes for an automorphism $\Phi \in \Out(G)$. We will see that like for free groups, there is a strong link between pseudo-periodic conjugacy classes and \emph{Nielsen paths}, which are some periodic paths in a train track representative for $\Phi$. Finding Nielsen paths in a train track representative will be our first goal. Then we establish the link with pseudo-periodic conjugacy classes.

\subsection{Computation of periodic indivisible Nielsen paths}

Like in the free group case, we can define Nielsen paths associated to a train track map. They are a useful tool to understand periodic elements of an automorphism. In the free group Nielsen paths are usually defined in the quotient graph but here we will define them in the tree instead, since in GBS graphs some non-degenerate turns cannot be seen in the quotient graph. Some example of references include \cite{BestvinaHandelTrainTracks}, \cite{BestvinaFeighnHandelLaminations}.

For any path $\alpha$ in a tree $T$ we define $[\alpha]$ as the path obtained by tightening $\alpha$ while keeping the same endpoints. Equivalently, it is the geodesic between the endpoints of $\alpha$. If the endpoint of $\alpha$ is the initial point of $\beta$ then $\alpha \cdot \beta$ is the concatenation of $\alpha$ and $\beta$. It is a \emph{tight} concatenation if $\alpha$ and $\beta$ are geodesics and $\alpha \cdot \beta$ is a geodesic. A tight concatenation is a \emph{legal} concatenation if the turn at the concatenation point is legal, otherwise it is \emph{illegal}.

Let $f: T \to T$ be an irreducible train track representative for $\Phi \in \Out(G)$. A \emph{Nielsen path} is a tight path $\gamma \subset T$ such that there exists $g \in G$ such that $f(\gamma)$ and $g \gamma$ are homotopic relative endpoints, i.e. $g^{-1} \cdot f$ fixes both endpoints of $\gamma$. 

A \emph{periodic Nielsen path}  is a tight path $\gamma \subset T$ such that there exists $n \geq 1$ and $g \in G$ such that $f^n(\gamma)$ and $g \gamma$ have the same endpoints. The minimal nonzero integer $n$ is the \emph{period} of $\gamma$.

\begin{rema}
 A tight concatenation of periodic Nielsen paths $\gamma_1 \cdot \gamma_2$ might not be a periodic Nielsen path but it is pre-periodic: there exists $k \in \N$ such that $[f^k(\gamma_1 \cdot \gamma_2)]$ is a periodic Nielsen path. This is a consequence of the local finiteness of trees: $[f^k(\gamma_1 \cdot \gamma_2)]$ has bounded length. Since the tree is locally finite and the orbits of the endpoints are periodic, there are finitely many possibilities for the orbit of $[f^k(\gamma_1 \cdot \gamma_2)]$, hence its periodicity.
\end{rema}

In general the endpoints of a Nielsen path are not vertices of the tree but rather lie in the interior of some edges.

An \emph{indivisible Nielsen path} (INP) is a Nielsen path which cannot be written as the tight concatenation of two shorter Nielsen paths. A \emph{periodic indivisible Nielsen path} (pINP) is a periodic Nielsen path which cannot be written as the concatenation of two shorter periodic Nielsen paths. Periodic INPs are defined as INPs for a certain power of $f$. Just as for periodic Nielsen path one can define the \emph{period} of a periodic INP.

Lemma \ref{lem:NielsenInclus} is a GBS version of well-known facts for pINPs in the free group case.

\begin{lem}\label{lem:NielsenInclus}\label{lem:formeINP}
 Let $f : T \to T$ be a train track map. 
 \begin{enumerate}[label=(\roman*)]
  \item Let $\gamma$ be a periodic Nielsen path of period $n$. There exists a unique decomposition of $\gamma$ as a tight concatenation of periodic indivisible Nielsen paths whose period divides $n$.
  \item If $\gamma$ is a tight concatenation of periodic indivisible Nielsen paths $\eta_1 \cdot \dots \cdot \eta_k$ and $\eta$ is a pINP such that $\eta \subset \gamma$ then there exists $i \in \{1, \dots, k \}$ such that $\eta=\eta_i$. In particular the decomposition of a periodic Nielsen path into pINPs is unique.
  \item A periodic indivisible Nielsen path contains a unique illegal turn. It is a tight concatenation $\alpha \cdot \beta$ where $\alpha, \beta$ are legal subpaths.
 \end{enumerate}
\end{lem}
\begin{proof}
 Let $\gamma$ be a periodic Nielsen path. Let $g \in G$ and $n \in \N$ such that $g \cdot f^n$ fixes the endpoints of $\gamma$.
 
 Let $\gamma_1, \dots, \gamma_k$ be the maximal legal subsegments such that $\gamma = \gamma_1 \cdot \dots \cdot \gamma_k$. The path $\gamma$ cannot be legal since the lengths of its images by $f^n$ is bounded so $k \geq 2$.
 
 We have $\gamma_i \subset g \cdot f^n(\gamma_i)$ for all $1 \leq i \leq k$. Indeed $f^n(\gamma_i)$ is a legal segment, so it can only intersect a unique $g^{-1} \gamma_j$, otherwise it would cross an illegal turn. Since $f^n(\gamma)$ covers $g^{-1} \gamma$ and $[f^n(\gamma)]$ has exactly the same number of illegal turns as $\gamma$, $f^n$ induces a permutation of the $\gamma_i$ for $1 \leq i \leq k$. By continuity of $f^n$ the permutation must be trivial.
 
 Up to taking a multiple of $n$ we may suppose that there is cancellation at every illegal turn when applying $f^n$. 
  For $2 \leq i \leq k-1$, $\gamma_i$ lies in the interior of $g \cdot f^n(\gamma_i)$. Since $f^n$ stretches legal segments uniformly this implies that there is a unique fixed point $v_i$ for $g \cdot f^n$ in the interior of $\gamma_i$. The vertices $v_i$ cut $\gamma$ into $k-1$ shorter periodic Nielsen paths, each of which has a unique illegal turn. 
 
 Since periodic indivisible Nielsen paths cannot be cut into smaller paths, they have at most one illegal turn, and at least one since their length does not grow exponentially with $f$. This proves (iii). Conversely, if a periodic Nielsen path has a unique illegal turn, it can be written as a concatenation of legal paths $\alpha \cdot \beta$. Up to replacing $f$ by $gf^n$, the first point of $\alpha$ and the last point of $\beta$ are fixed by $f$ and all other points in $\alpha \cdot \beta$ escape exponentially when iterating $f$. Therefore pINPs are exactly periodic Nielsen paths with one illegal turn.
 For a general periodic Nielsen path $\gamma$, we found a decomposition into periodic Nielsen paths with one illegal turn each, also pINPs so (i) is proved.
 
 \bigskip
 
 Let us prove (ii). Now we do not assume that $\gamma$ is a periodic Nielsen path, only that it is a concatenation of pINPs. We just proved that the former implies the latter. Observe that if $\alpha, \beta$ are legal paths and if $\eta = [x, y]$ is a pINP with $x \in \alpha$, $y \in \beta$ then for any $a, b \in \alpha \cdot \beta$ such that $\{a, b\} \neq \{x, y\}$, either $d(f^n(a), f^n(b)) \rightarrow \infty$ or there exists $n \in \N$ such that $f^n(a)=f^n(b)$. 
 
 Let $\eta = [x,y]$ be a pINP in $\gamma$. 
 There exists a unique illegal turn in $\eta$ and it is also the illegal turn in a pINP $\gamma'= [x', y']$ which appears in a decomposition of $\gamma$.
 
 Let $\alpha, \beta$ be the maximal legal subsegments in which $x, y$ lie. Since $\gamma'$ has a unique illegal turn, $x',y'$ lie in $\alpha \cdot \beta$. Since $\gamma'$ is a pINP the sequence $d(f^n(x'), f^n(y'))$ is bounded and positive. The claim implies that $\{x',y'\} =\{x,y\}$.
 
 This proves the uniqueness of the decomposition.
\end{proof}

Recall the bounded cancellation constant $\BCC(f)$. Suppose $\alpha \cdot \beta$ is the concatenation of two legal paths.
Then we have $\len([f(\alpha \cdot \beta)]) > \lambda \len(\alpha \cdot \beta) - 2 \BCC(f)$.

We can deduce a bound on the length of Nielsen paths.  Let $C_f := \frac{2 \BCC(f)}{\lambda - 1}$.  If $\alpha, \beta$ are legal and $\len (\alpha \cdot \beta) > C_f$ then $\len ([f(\alpha \cdot \beta)]) > \len (\alpha \cdot \beta)$ so the length of $([f^n(\alpha \cdot \beta)])$ is strictly increasing. Therefore, if $\alpha \cdot \beta$ is an indivisible Nielsen path, then $\len (\alpha) \leq C_f, \len(\beta) \leq C_f$.

The following fact is a well-known fact for the free group case and also applies here, mainly because trees are also locally finite in the GBS case.
\begin{lem} \label{lem:finiNielsen}
 There are finitely many orbits of periodic indivisible Nielsen paths.
\end{lem}

\begin{rema}
 Since $T$ is locally finite, there are finitely many orbits of edge paths with bounded length. However in general the endpoints of pINPs are not vertices of $T$ so pINPs are not edge paths.
\end{rema}

 Before proving Lemma \ref{lem:finiNielsen} we introduce the following notion which will enable us to work with edge paths while looking for periodic indivisible Nielsen paths.
 \begin{defi}
  A \emph{pseudo-pINP} is an edge path $\gamma \subset T$ with one illegal turn such that there exists $n \geq 1$ and $g \in G$ such that $\gamma \subset g[f^n(\gamma)]$.
 \end{defi}
 \begin{lem} \label{lem:correspondance-chemin-arete}
  Every pseudo-pINP contains a unique pINP. Conversely, the minimal edge path containing a pINP is a pseudo-pINP.
 \end{lem}
 \begin{proof}
  Suppose $\gamma$ is an edge path with one illegal turn, and let $n \geq 1$ and $g \in G$ be such that $\gamma \subset g[f^n(\gamma)]$. Up to replacing $f$ by $gf^n$ we may suppose $\gamma \subset [f(\gamma)]$. Write $\gamma$ as the concatenation $\alpha \cdot \eta_1 \cdot \eta_2 \cdot \beta$ where the illegal turn is between $\eta_1$ and $\eta_2$, and $\eta_1$ and $\eta_2$ are the maximal subpaths such that $f(\eta_1)= f(\bar \eta_2)$. If $\gamma$ contains a pINP $\gamma'$ then the length of each legal branch of $\gamma'$ must be $\frac{\lambda \len(\eta_1)}{\lambda-1}$, which proves uniqueness. For the existence, since $\gamma \subset f(\gamma)$ this length defines a unique subpath $\gamma' \subset \gamma$ and $[f(\gamma')]=\gamma'$. 
  
  Conversely, the minimal edge path $\gamma$ containing a pINP $\gamma'$ is sent to an edge path containing $gf^n(\gamma')$ for some $n \geq 1$ and $g \in G$. Since $\gamma' \subset g[f^n(\gamma')] \subset g[f^n(\gamma)]$ and $\gamma$ is the minimal edge path containing $\gamma'$, then $\gamma$ is a pseudo-pINP.
 \end{proof}

 Now we prove Lemma \ref{lem:finiNielsen}.
\begin{proof}[Proof of Lemma \ref{lem:finiNielsen}]
 First of all the tree $T$ is locally finite and there are finitely many orbits of illegal turns. There is a bound on the lengths of pINPs, namely the critical constant $C_f$. Finally every pINP crosses a unique illegal turn.

 Every pINP is contained in a pseudo-pINP whose legal branches have length at most $C_f + \max_{e \in E(T)} \len(e)$. There are finitely many orbits of edge paths with bounded length, hence finitely many orbits of pseudo-pINPs with bounded length. Since every pseudo-pINP contains a unique pINP by Lemma \ref{lem:correspondance-chemin-arete}, there are finitely many orbits of pINPs.
\end{proof}

 The important implication of Lemma \ref{lem:finiNielsen} is that the set of orbits of indivisible Nielsen paths can be computed algorithmically. To prove this for GBS groups we use an approach which resembles \cite[Sections 5, 6]{CoulboisLustigLongTurns}, where \emph{long turns} are used to understand pINPs. 
 
 The algorithm which finds all periodic indivisible Nielsen paths relies on the correspondance with pseudo-pINPs given by Lemma \ref{lem:correspondance-chemin-arete}. It suffices to prove that all pseudo-pINPs can be computed algorithmically. This would be quite simple if we knew the bound $C_f$. However we do not know how to estimate $\BCC(f)$ and instead we will give an algorithm which does not need it. The algorithm consists in testing all edge paths with bounded length for pseudo-pINPs, and using a criterion to make sure that the bound on the length was sufficiently big. If not then the bound is increased: eventually it becomes greater than $C_f$ and the criterion is satisfied.

 \bigskip
  
 Let $\gamma = \alpha \cdot \beta$ be an edge path with $\len(\alpha), \len(\beta) \neq 0$ such that $\alpha, \beta$ are legal subpaths and the turn at the concatenation is illegal. 
 
 The following technical lemma studies the behaviour of the sequences  $([f^n(\gamma)])_{n \in \N}$.
 
\begin{figure}
  \centering
  \includegraphics[scale=0.4]{\localhost/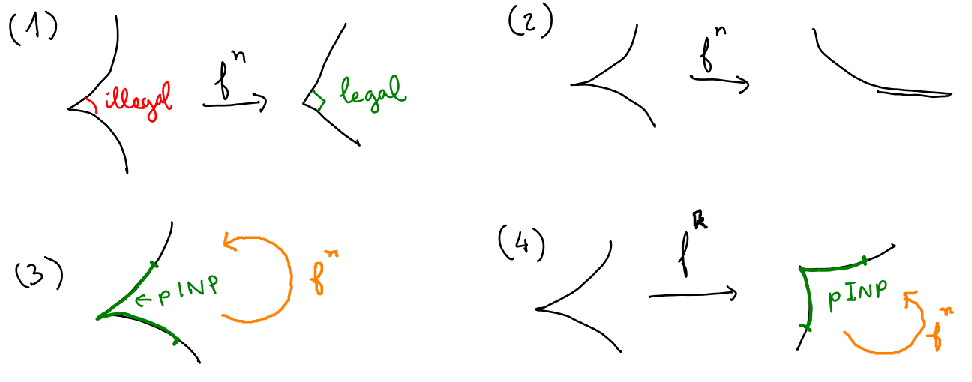}
  \caption{Four different cases for the behaviour of $[f^n(\gamma)]$} \label{fig:4comportements}
 \end{figure}

 \begin{lem}\label{lem:techniqueSuiteGamma}
 \begin{enumerate}[label=(\roman*)]
  \item There are four mutually exclusive behaviours for the sequence $([f^n(\gamma)])$, illustrated by Figure \ref{fig:4comportements}:
      \begin{enumerate}[label=(\arabic*)]
       \item for some $n \in \N$, $[f^n(\gamma)]$ is legal and neither $f^n(\alpha)$ nor $f^n(\beta)$ contains the other
       \item for some $n \in \N$, $f^n(\alpha) \subset f^n(\beta)$ or $f^n(\beta) \subset f^n(\alpha)$
       \item $\gamma$ contains a pseudo-pINP
       \item $\gamma$ does not contain any pseudo-pINP but there exists $k \geq 1$ such that $[f^k(\gamma)]$ does
      \end{enumerate}
  \item There is an algorithm which takes an edge path $\gamma$ with one illegal turn and returns what case $([f^n(\gamma)])$ belongs to
  \item If $\alpha, \beta$ are both longer than $C_f$ then case (2) is not possible
  \item If there exists a pINP $\eta$ such that the minimal edge path containing $\eta$ strictly contains $\gamma$ then $\gamma$ is in case (2).
 \end{enumerate}
 \end{lem}
 
 \begin{proof}
  (i) If there exists $n \in \N$ such that $[f^n(\gamma)]$ is legal then the sequence belongs to case (1) or (2). Otherwise, for all $n \in \N$, the path $[f^n(\gamma)]$ contains an illegal turn. In particular its length is never zero. Let $B \geq C_f$ be a bound on the length of legal branches of pseudo-pINPs. There are finitely many edge paths with length at most $2B$ so there exist $k\geq 0, n \geq 1$ and $g \in G$ such that $[f^k(\gamma)]$ and $g[f^{k+n}(\gamma)]$ are either equal, or intersect on a length at least $B$ on each side of the illegal turn. Let $\eta$ be the maximal edge path in the intersection $[f^k(\gamma)]\cap g[f^{k+n}(\gamma)]$. Let $\alpha_\eta$ be one of the legal branches of $\eta$. Since $\len(\alpha_\eta) \geq C_f$ we have $\len(f^n(\alpha_\eta) \cap [f^n(\eta)]) \geq \len(\alpha_\eta)$. Therefore the image $[f^n(\eta)]$ contains $g^{-1} \eta$, so $\eta$ is a pseudo-pINP. The sequence belongs to case (3) if $k = 0$ and to case (4) otherwise.
  
  Case (1) obviously excludes cases (2), (3) and (4). Case (2) excludes cases (3) and (4) since it implies that $[f^n(\alpha \cdot \beta)]$ is a legal path. Definition of cases (3) and (4) also implies that they are mutually exclusive.

 (ii) Before giving the algorithm which determines the case, let us give some preliminaries:
 \begin{itemize}
  \item By Proposition \ref{prop:algos-de-base} (v), there is an algorithm which, given finite non-backtracking edge paths $\alpha, \beta$, decides whether they belong to the same orbit or whether there exists $g \in G$ such that $\alpha \subset g \beta$. 
  \item There is an algorithm which takes an edge path $\eta$ and a period $p \geq 1$, and determines if $\eta$ is a pseudo-pINP with period at most $p$. Compute the iterates $[f^n(\eta)]$ for $n \in \{1, \dots, p\}$. If there exists $n$ such that $[f^n(\eta)]$ contains a translate of $\eta$ then $\eta$ is a pseudo-pINP whose period divides $n$. Otherwise it is not a pseudo-pINP with period at most $p$.
 \end{itemize}

   Let $\gamma$ be an edge path with one illegal turn. Repeat the following steps, starting with $n= 0$:
    \begin{itemize}
     \item Compute $f^{n}(\gamma)$.
     \item If $f^{n}(\gamma)$ is legal, check whether $f^{n}(\alpha) \subset f^{n}(\beta)$ or vice-versa. If yes then $\gamma$ belongs to case (2), else it belongs to case (1). The algorithm stops.
     \item If $f^{n}(\gamma)$ is not legal, check if it contains a pseudo-pINP with period at most $n+1$. Suppose it contains a pseudo-pINP with period $p \leq n+1$. Then either case (3) or case (4) is true. The pINP contained in $\eta$ also has period $p$: if $\gamma$ contains a pINP $\eta_0$ then $\eta$ contains its $f^n$-image, so $\eta_0$ has period $p$. Then it suffices to check if $\gamma$ contains any pseudo-pINP with period $p$. If yes then $\gamma$ is in case (3), else it is in case (4).
     \item Replace $n$ by $n+1$.
    \end{itemize}
  Because of (i) this will eventually terminate.
 
 (iii) If $\len(\alpha)\geq C_f$ then for every $n \in \N$, $\len(f^n(\alpha) \cap [f^n(\gamma)]) >\lambda  \len(\alpha) -\BCC(f)> \len(\alpha)$. A similar statement holds if $\len(\beta) \geq C_f$. When both $\alpha$ and $\beta$ are longer than $C_f$, (2) cannot occur.
  
  (iv) Suppose that there exists a pINP $\eta$, a pseudo-pINP $\eta'$ containing $\eta$, such that $\gamma \subsetneq \eta'$, and that $\gamma$ does not contain $\eta$. One of the endpoints $x$ of $\gamma$ lies in the interior of $\eta$. Let $y$ be the endpoint of $\eta$ in the same legal branch as $x$. The distance $d_T(f^n(x), f^n(y))$ increases exponentially since $[x, y]$ is legal. The distance between $f^n(y)$ and the illegal turn is bounded by $C_f$. For $n_0$ big enough, $d_T(f^{n_0}(x), f^{n_0}(y))> C_f$ so $f^{n_0}(x)$ lies in the simplified part of $f^{n_0}(\eta)$. Since $\gamma \subset \eta'$, $[f^{n_0}(\gamma)]$ is contained in $[f^{n_0}(\alpha)]$ or $[f^{n_0}(\beta)]$ so $\gamma$ is in case (2).
 \end{proof}

\begin{prop} \label{prop:compute-pinps}
 There is an algorithm taking a train track map $f : T \to T$ and finding all orbits of minimal pseudo-pINPs.
\end{prop}
\begin{proof}
 Proposition \ref{prop:algos-de-base} (v) implies that given two edge paths, one can decide whether they belong to the same orbit. As a result one can list all orbits of edge paths of given length: it suffices to choose a representative for each orbit of vertex, and then construct all edge paths with given length starting with these vertices.
 
 Start with $L = 1$. Apply the following steps.
 \begin{itemize}
  \item Let $\L_L$ be a list of representatives of edge paths with shape $\alpha \cdot \beta$ where $\alpha$ and $\beta$ are legal paths with combinatorial length at most $L$.
  \item For each path in $\L_L$, determine what case they belong to using Lemma \ref{lem:techniqueSuiteGamma} (ii).
  \item If there exist a path in case (2) then increase $L$ by $1$ and start again. Else stop. Let $\mathcal T$ be the subset of $\L_L$ consisting of paths in case (3). 
 \end{itemize}
 Every path $\gamma \in \mathcal T$ contains a pseudo-pINP, moreover point (iii) of Lemma \ref{lem:techniqueSuiteGamma} ensures that we can find such a pseudo-pINP $\eta_0 \subset \gamma$ and its period $p$. For every other pseudo-pINP $\eta \subset \gamma$, the period of $\eta$ is also $p$. Compute $[f^p(\eta)]$ for every $\eta \subset \gamma$ in order to find all pseudo-pINPs in $\gamma$ and find the minimal one.
 
 \bigskip
 
 This algorithm eventually stops because after enough steps $L$ becomes greater than $C_f$ so case (2) does not occur. Suppose that the algorithm stops, then for every minimal pseudo-pINP $\eta$ there exists $\gamma \in \mathcal T$ such that $\eta \subset \gamma$. By contradiction suppose otherwise. Then $\len(\eta) > L$. There exists a subpath $\gamma \subsetneq \eta$ of combinatorial length $L$ and by minimality of $\eta$, $\gamma$ does not contain the pINP contained in $\eta$. Then by (iv) of Lemma \ref{lem:techniqueSuiteGamma} $\gamma$ is in case (2) which is a contradiction.
\end{proof}

Of course it is more convenient to work with edge paths than with arbitrary paths which can start or end in the middle of an edge, especially for algorithmic purposes. Instead of using pseudo-pINPs, we will subdivide $T$ at the endpoints of pINPs so that all pINPs become actual edge paths:
\begin{lem} \label{lem:existe-subdivision}
 There exists a subdivision $s:T \to T'$ with $T' \in \D$  and a train track map $f' : T' \to T'$ such that $f' \circ s = s \circ f$, such that the endpoints of all pINPs for $f'$ are vertices of $T'$.
\end{lem}
\begin{proof}
 Define $T'$ by subdividing $T$ at every endpoint of periodic indivisible Nielsen paths. Since there are finitely many orbits of pINPs the tree $T'$ is simplicial. 
 Let $s: T \to T'$ be the corresponding isometry.
 
 Since the set of pINPs is stable under $f$, the map $f'$ induced by $f$ on $T'$ maps vertex to vertex. The map $f'$ is a train track map.
\end{proof}

\begin{lem} \label{lem:compute-subdivision}
 One can compute a subdivision $s : T \to T'$ and a train track map $f': T' \to T'$ such that $f' \circ s = s \circ f$ and such that all pINPs for $f'$ start and end at vertices of $T'$.
\end{lem}
\begin{proof}
 The map $f$ exists by Lemma \ref{lem:existe-subdivision}. We will give an algorithmic construction.
 
 One can find all orbits of minimal pseudo-pINPs by Proposition \ref{prop:compute-pinps}. Let $\gamma$ be a minimal pseudo-pINP. It contains a unique pINP $\eta$. There exists $n \geq 1$ and $g \in G$ such that $\eta = g[f^n(\eta)]$. The first point of $\eta$ (resp. last point) is a vertex if and only if the first edge (resp. last edge) of $\gamma$  is equal to the first edge (resp. last edge) of $g[f^n(\gamma)]$. This can be checked algorithmically.
 
 There exists $k \geq 1$ such that all pINPs are $k$-periodic. Step by step, we will construct a subdivision of $T$ which is $f^k$-invariant. We will check that it is also $f$-invariant.

 Start with $S = T$.
 
 While there exists a pINP $\eta=[x,y]$ for $f^k$ such that $x$ is not a vertex, define the subdivided tree $S'$ as follows. Let $e$ be the edge of $S$ which contains $x$. Subdivide $e$ by adding a new vertex $v$ representing $x$. The map $f^k$ induces a map on $S'$ with $f^k(x) = gx$ where $g$ is such that $f^k(\eta)=g \eta$.
 
 By repeating this with $S:=S'$ for every pINP, we obtain an $f^k$-invariant subdivision $T'$.
  
 \bigskip
 
 Now we would like to prove that the subdivision is $f$-invariant. We need to define $f$ on vertices of $V(T') \setminus V(T)$.  Suppose $v \in V(T') \setminus V(T)$. In $T$ the point $v$ is not a vertex but it is the endpoint of a pINP $\eta$. There exists a unique pINP $\eta'$ such that $[f(\eta)] = \eta'$, and the path $\eta'$ can be computed. The map $f$ sends the endpoints of $\eta$ to those of $\eta'$: define $f'(v)$ as the corresponding endpoint of $\eta'$ in $T'$. By construction $f'(v)$ is a vertex.

 The construction does not depend on the choice of $\eta$ since it exhibits the subdivision of Lemma \ref{lem:existe-subdivision}.
\end{proof}

 \subsection{Pseudo-periodic conjugacy classes of $G$}

 Recall that the conjugacy class of $g \in G$ is \emph{pseudo-periodic} for $\Phi \in \Out(G)$ if $\| \phi^n(g)\|_T$ is bounded for $n \rightarrow \infty$ and $\phi \in \Phi$.
 
  The motivation for the notion of pseudo-periodicity is given by the following result, which we will prove below:
 \begin{prop} \label{prop:lien-classe-conj-facteur}
  If there exists a simple element $h \in G$ such that $h$ is pseudo-periodic, then a power of $\Phi$ is reducible.
 \end{prop}
 
 \begin{lem}\label{lem:nonInfiniBornee}
  Let $g \in G$ be a loxodromic element. Let $\phi \in \Aut(G)$:
  The following conditions are equivalent:
    \begin{itemize}
     \item $\| \phi^n(g)\|_T$ does not tend to infinity
     \item $\| \phi^n(g)\|_T$ is bounded.
    \end{itemize}
 \end{lem}
 \begin{proof}
  The axis of $g$ in $T$ can be partitioned into maximal legal segments, concatenated at illegal turns. Let $s_1, \dots, s_m$ be consecutive maximal legal segments forming a fundamental domain of $\axe_T(g)$.
  
  For every $n \in \N$ there exists $I_n \subset \{1, \dots, m\}$ such that $\bigcup_{i \in I_n} f^n(s_i)$ contains a fundamental domain of the axis of $\phi^n(g)$. We may assume $I_n \subset I_{n-1}$ for every $n \in \N$. For $i \in I_n$ let $J_i^n := f^n(s_i) \cap \axe_T(\phi^n(g))$. It is a legal segment.

  If $\|\phi^n(g)\|_T$ is not bounded then at least one of the $J_i^n$ has unbounded length. Suppose for $N \in \N$, $\len(J_{i_0} ^N) > 2C_f$. When applying $f$ to $J_{i_0}^n$, it is stretched by a factor $\lambda$, however there might be cancellation at the ends of $J_{i_0}^n$ because of illegal turns. Since this cancellation cannot exceed $\BCC(f)$, we have $\len(J_{i_0}^{N+1}) \geq \lambda \len(J_{i_0}^N)  - 2 \BCC(f) > \frac{\lambda}{2}  \len(J_{i_0}^N) +  \lambda  C_f - 2\BCC(f) > \frac{\lambda}{2}  \len(J_{i_0}^N)$ so $\len(J_{i_0}^n) \geq  \left ( \frac{\lambda}{2} \right ) ^{n-N}  \len(J_{i_0}^N)$. Therefore $\| \phi^n(g)\|_T$ goes to infinity.
 \end{proof}

 \begin{rema}
  Equivalently the conjugacy class of an element $g \in G$ is pseudo-periodic for $\Phi$ if there exists $S \in \D$ such that $\| \phi^n(g)\|_S$ is bounded with $\phi \in \Phi$, and equivalently if for all $S \in G$, $\| \phi^n(g)\|_S$ is bounded.
 \end{rema}

 Pseudo-periodic elements are the right analogue of periodic elements in the GBS context. 
  Although the conjugacy class of a periodic element is not periodic, the translation length is.
 \begin{lem}\label{lem:2domainesFonda}
 Let $\phi \in \Aut(G)$. 
 Let $h,g \in G$ be loxodromic elements. Suppose there exists $t \in \D$ and $v \in \axe_T(h) \cap \axe_T(g)$ such that $[v,h^2v]=[v,g^2v]$. Then in any $S \in \D$ there exists $w \in S$ such that $[w, g^2w] = [w, h^2w]$ is a common pair of fundamental domains for $g$ and $h$. In particular $\|g\|_S = \|h\|_S$ for all $S \in \D$.
\end{lem}

\begin{figure}
 \centering
   \includegraphics[scale=0.35]{\localhost/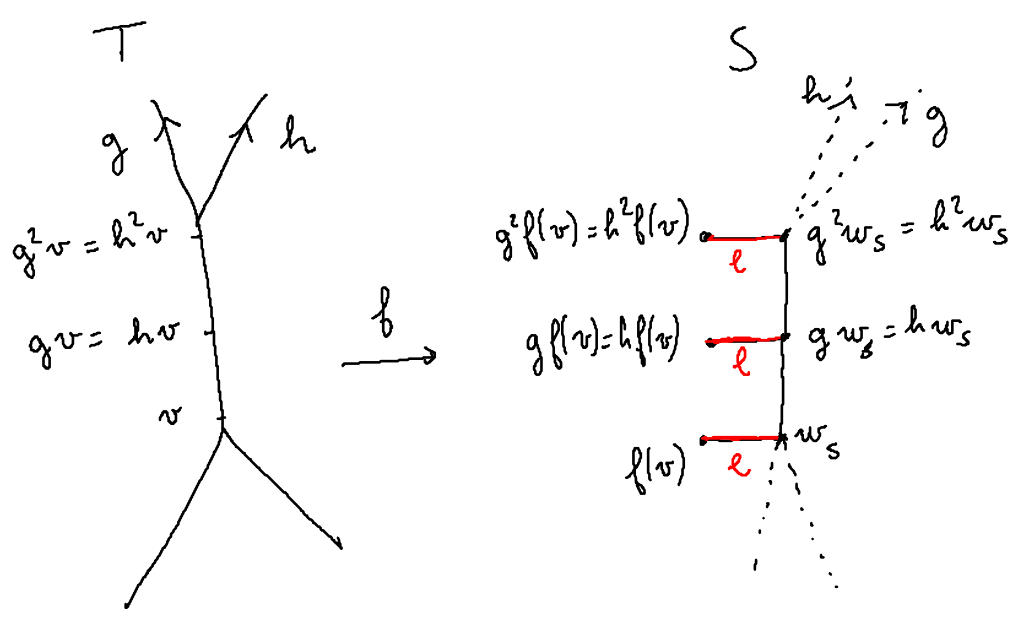}
 \caption{Description of the proof of Lemma \ref{lem:2domainesFonda}: image of the pair of fundamental domains in the new tree $S$} \label{fig:deuxDomainesFonda}
\end{figure}

\begin{proof}
 Let $\phi, h,g, T ,v$ be as described in the hypotheses. 

 Let $f:T \to S $ be a $G$-equivariant application. See Figure \ref{fig:deuxDomainesFonda}. The segment $[f(v), hf(v)]$ intersects $\axe_S(h)$. Let $l_h$ be the distance between $f(v)$ and $\axe_S(h)$. It is equal to the distance between $hf(v)$ and the axis. Moreover the length of the intersection $[f(v),hf(v)] \cap [hf(v), h^2f(v)]$ is equal to $l_h$. The same goes for $l_g$ and since $gv = hv$ and $g^2v = h^2v$, we have $l_g = l_h$.
 
 Let $w_S \in S$ be the point of $[f(v), hf(v)]$ at distance $l_h=l_g$ from $f(v)$. It belongs to $\axe_S(h)$ and $\axe_S(g)$. 
 
 The translation length of both $g$ and $h$ in $S$ is $d_S(f(v), hf(v)) - 2 l_h$. Therefore $hw_S = gw_S$ is the point of $[f(v), hf(v)]$ at distance $l_h$ from $hf(v)$ and $h^2 w_S= g^2w_S$ is the point at distance $l_h$ from $h^2(v)$. 
\end{proof}

\begin{rema} \label{rem:N-domaines-fonda}
 The previous lemma is also true when replacing $2$ with any $n \geq 2$: if $g, h$ share $n$ consecutive fundamental domain in some $T \in \D$ then they do in every $T \in \D$.
\end{rema}

\begin{cor} \label{coro:2domainesFonda}
 Let $h,g \in G$ be loxodromic elements such that there exists $T \in \D$ and $v \in \axe_T(h) \cap \axe_T(g)$ such that $[v,h^2v]=[v,g^2v]$. Then any special factor containing $h$ also contains $g$.
\end{cor}

\begin{proof}
 Let $A$ be a special factor containing $h$. Let $S \in \D^\A$ be a tree such that a collapse $\bar S$ of $S$ has a vertex $x$ with stabilizer $A$. By Lemma \ref{lem:2domainesFonda} the axes of $h$ and $g$ in $S$ share a fundamental domain $[w, gw] \subset S$.
 
 In $\bar S$, $A$ is the stabilizer of the vertex $x$. Therefore the axis of $h$ is sent to $x$ by the collapse map $S \to \bar S$. The fundamental domain of $g$ in the axis of $h$ is also sent to $x$, so the whole axis of $g$ is sent to $x$ by equivariance. Thus $g$ fixes $x$ and $g \in A$.
\end{proof}

In the following result, the fact that trees in $\D$ are locally finite is fundamental.

\begin{prop} \label{prop:equivalencePseudoperiodique}
 Let $\phi \in \Aut(G)$. Let $g \in G$ be a loxodromic element. The following conditions are equivalent:
 \begin{enumerate}[label=(\roman*)]
  \item the conjugacy class of $g$ is pseudo-periodic for the outer class $\Phi = [\phi]$
  \item for any $T \in \D$, for all $N \geq 1$, there exist distinct $n, m \in \N$ such that $\phi^n(g)$ and $\phi^m(g)$ admit $N$ common consecutive fundamental domains up to translation
  \item for all $N \geq 1$, there exists $m>0$  such that for any $S \in \D$,  $g$ and $\phi^m(g)$ admit $N$ common consecutive fundamental domains up to translation in $S$
  \item for any $T \in \D$ the sequence $\left ( \| \phi^n(g)\|_T \right )_{n \rightarrow - \infty}$ is bounded.
 \end{enumerate}
\end{prop}
\begin{proof}
 Let us prove (i) $\Rightarrow$ (ii). Suppose the conjugacy class of $g$ is pseudo-periodic.
 
 Let $f: T \to T$ be a representative for $\phi$. Let $B > \sup_{n \in \N} \| \phi^n(g)\|_T$. Let $N\in \N$. 
Since there are finitely many orbits of edge paths with length at most $NB$, there exist distinct $n,m \in \N$ and $h \in G$ such that $\phi^n(g)$ and $h \phi^m(g) h^{-1}$ share $N$ consecutive fundamental domains.
 
 To prove (ii) $\Rightarrow$ (iii), take $T \in \D$. By (ii) there exists $n,m$ such that $\phi^n(g)$ and $\phi^m(g)$ have $N$ consecutive fundamental domains up to translation. By Lemma \ref{lem:2domainesFonda} and Remark \ref{rem:N-domaines-fonda}, these two elements have $N$ consecutive fundamental domains, up to translation, in every tree in $S \in \D$. In particular in $S \cdot \phi^{-n}$ they do. This implies that $g$ and $\phi^{m-n}(g)$ share $N$ fundamental domains up to translation in $S$, for every $S \in \D$.

 Suppose (iii): in particular, by taking $S' := S \cdot \phi^n$ we have $\|\phi^n(g)\|_S = \| \phi^{n+m}(g)\|_S$ for any $n \in \Z$. Therefore $\| \phi^n(g)\|_S$ is bounded so $g$ is pseudo-periodic for $[\phi]$ so (i) and (iv) are true.
 
 Finally suppose (iv). Condition (iv) is condition (i) for $\phi^{-1}$ instead of $\phi$, and it implies (iv) for $\phi^{-1}$, which is (i) for $\phi$. Therefore (iv) $\Rightarrow$ (i).
\end{proof}

Now we can prove Proposition \ref{prop:lien-classe-conj-facteur}, which states that the minimal factor containing an element whose conjugacy class is pseudo-periodic is itself periodic.

 \begin{proof}[Proof of Proposition \ref{prop:lien-classe-conj-facteur}]
  Condition (iii) of Proposition \ref{prop:equivalencePseudoperiodique} states that there exist $n,m$ such that $h$ and $\phi^m(h)$ share two consecutive fundamental domains up to translation. Using Corollary \ref{coro:2domainesFonda} this implies that any special factor containing $h$ also contains a conjugate of $\phi^m(h)$, and vice versa. Thus the minimal factors 
  containing these elements are conjugate and the conjugacy class of these factors is $\phi^{m}$-periodic.
 \end{proof}

 \bigskip

 The following results are the key for finding pseudo-periodic elements in $G$: 
 \begin{prop} \label{prop:concatenationNielsen}
  Let $g \in G$ be an element whose conjugacy class is pseudo-periodic for $\Phi \in \Out(G)$. Suppose $f: T \to T$ is a train-track representative for $\Phi$. Then the axis of $g$ in $T$ is a concatenation of periodic indivisible Nielsen paths.
 \end{prop}
 \begin{proof}
   Let $\phi \in \Phi$ be such that $f$ is $\phi$-equivariant.
   The axis of $g$ in $T$ is a concatenation of maximal legal subsegments interrupted by illegal turns. The tightened image of the axis of $g$ by $f^n$ is the axis of $\phi^n(g)$. Since $f^n$ maps legal segments to legal segments, the number $N_n$ of orbits of maximal legal segments in $\axe_T(\phi^n(g))$ never increases with $n$. Since $g$ and $\phi^{nk}(g)$ share two fundamental domains up to translation (Proposition \ref{prop:equivalencePseudoperiodique} (iii)) for some fixed $k$ and every $n \in \N$, the number $N_n$ is actually a constant $N$.

   By Proposition \ref{prop:equivalencePseudoperiodique} (iii) there exists $n \geq 1$ and $h \in G$ such that $h \phi^n(g) h^{-1}$ and $g$ have at least $N+1$ consecutive fundamental domains in common.   Up to replacing $\phi$ by $c_{h } \circ \phi^n$ for some $l \in \Z$, and $f$ by $ h f^n$, we may suppose that $g$ and $\phi(g)$ share $N+1$ consecutive fundamental domains. Let us call $\sigma$ the segment where both axes overlap. Up to replacing $\phi$ by $c_{\phi(g^l)} \circ \phi $ we may suppose that for any fundamental domain $\eta \subset \sigma$, the first point of $f(\eta)\cap \sigma$ is contained in $\eta$ (see Figure \ref{fig:segmentsEtAxe}). Note that the set of pINPs for $\phi$ does not change when replacing $\phi$ with a power or composing with an inner automorphism, so these assumptions will not change the outcome of the proof.
  
   \begin{figure}
    \centering
    \includegraphics[scale=0.4]{\localhost/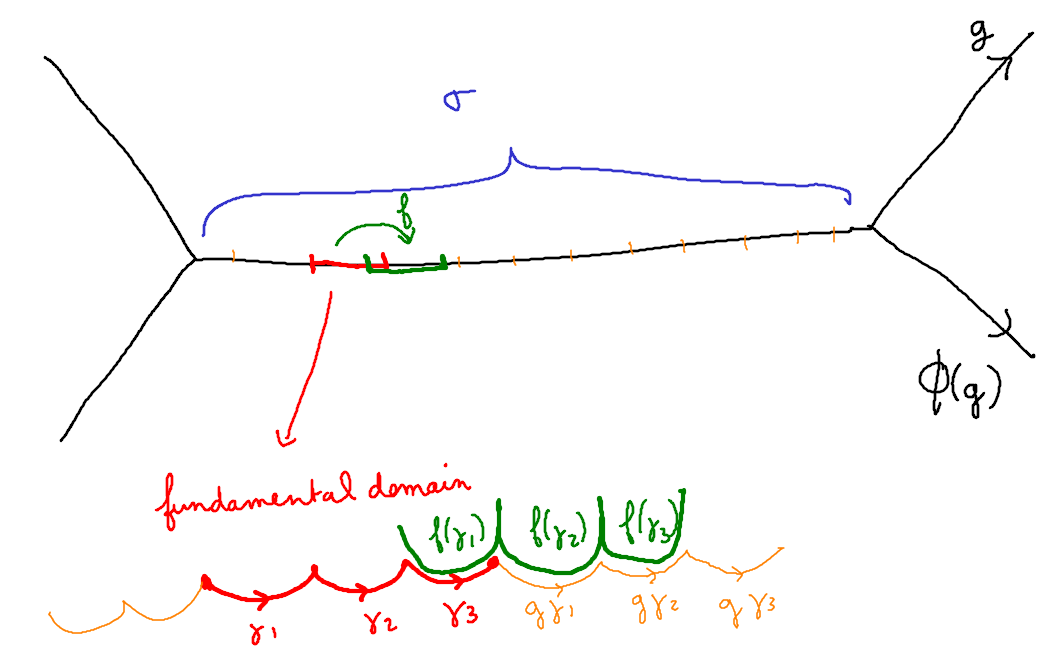}
    \caption{Picture of the action of $f$ on legal segments} \label{fig:segmentsEtAxe}
   \end{figure}
   
   Let $\gamma_1, \dots, \gamma_N, \dots, \gamma_{N^2}, \gamma_{N^2 + 1}$ be consecutive maximal legal segments which appear whole in $\sigma$. The isometry $g$ shifts the legal subsegments: for $i \in \{1, \dots, N^2 + 1 - N\}$ we have $\gamma_{i+N} = g \gamma_i = \tilde \phi(g) \gamma_i$.
   
   The map $f$ sends legal paths to legal paths, and it induces a bijection between the sets of maximal legal subsegments of $\axe_T(g)$ and $\axe_T(\phi(g))$. This bijection preserves the order of the segments. On the part where the axes overlap, $f$ shifts the subsegments $\gamma_i$ by an amount which does not depend on $i$. More accurately, there exists $0 \leq j < N$ such that for every $1 \leq i \leq N(N-1)$, $f(\gamma_i) \supset \gamma_{i + j}$.
   
   In particular, for $i \in \{ 1, \dots, N + 1\}$ we have 
   \[
    f ^N(\gamma_i) \supset f^{N-1}(\gamma_{i+j}) \supset \dots \supset \gamma_{i+ Nj} =g^j\gamma_i.
   \] 
   
   Since for all $i \in \{1, \dots, N +1\}$ we have $g^j \gamma_i \subset f ^N(\gamma_i)$ and $f^N$ stretches legal segments uniformly, there exists a unique point $p_i \in \gamma_i$ such that $f^N(p_i)= g^j p_i$. Besides $p_{N+1} = g p_1$. Thus for $1 \leq i \leq N$ the segment $[p_i, p_{i+1}]$ is a periodic Nielsen path, and since it has a unique illegal turn it is a pINP for $f$. Therefore $\axe_T(g)$ is a concatenation of pINPs.
 \end{proof}

 We also have a converse:
 \begin{prop}
  Suppose that $g \in G$ is a loxodromic element whose axis is a concatenation of pINPs. Then $g$ is pseudo-periodic.
 \end{prop}
 \begin{proof}
  Suppose the axis of $g \in G$ is a concatenation of pINPS. There exists a fundamental $\eta_1 \cdot \dots \cdot \eta_k$ for $g$ where $\eta_i$ is a pINP for $i \in \{1, \dots, k\}$. There exists a common $n \in \N$ such that for all $i \in \{1, \dots, k\}$ there exists $h_i \in G$ such that $[f^n(\eta_i)] = h_i \eta_i$. By continuity of $f^n$ the paths $h_i \eta_i$ and $h_{i+1} \eta_i$ have a common endpoint, although $h_i$ may be different from $h_{i+1}$. 
  
  Therefore, for every $k \in \N$, $\len([f^{nk}(\eta_1 \cdot \dots \cdot \eta_k)]] \leq \len(\eta_1 \cdot \dots \cdot \eta_k)$. Since $f$ is Lipschitz, this proves that $\|\phi^m(g)\|_T \rightarrow \infty$ when $n \rightarrow \infty$.
 \end{proof}
 \begin{rema}
  A consequence is that there cannot be any illegal concatenation of pINPs forming the axis of a loxodromic element of $G$. If there were such an element $g$, up to taking a power of $\Phi$ we could suppose the axis of $g$ is made of INPs with an illegal concatenation occuring at some point. Then the sequence $\| \phi^n(g) \|_T$ would be non increasing and non constant because of the illegal turn. This contradicts Proposition \ref{prop:equivalencePseudoperiodique} (iv).
 \end{rema}

\section{Pseudo-periodic subgroups} \label{sec:groupes-pseudo}
 
 In this section we introduce a collection of subgroups whose loxodromic elements are exactly the elements of $G$ whose conjugacy class is pseudo-periodic, and we give an algorithm which computes these subgroups.
 
 These groups can be understood as \emph{pseudo-periodic subgroups}, which are analogous to fixed subgroups for free groups. The notion of pseudo-periodicity which is used here depends on a choice of actual automorphism $\phi \in Aut(G)$ and is not a conjugacy class invariant.

 \subsection{Algorithmic computation of pseudo-periodic conjugacy classes}
 
 Fix $\Phi \in \Out(G)$.  Our aim in this subsection is to determine a subset of points which can be effectively constructed, and such that the loxodromic elements of its stabilizer are pseudo-periodic. 

Let $f: T \to T$ be a train track map for $\Phi$. Let $\eta :=[x_0, x_1] \subset T$ be a periodic indivisible Nielsen path for $f$. The \emph{Nielsen class} associated to $\eta$ is the set
 \[
  VY(\eta) := \{y \in T / \exists x_0 = y_0, y_1, \dots, y_n \text{ s.t. } \forall 1 \leq i \leq n \quad [y_{i-1}, y_i] \text{ is a pINP } \}
 \]

 The set of pINPs with both ends in $VY(\eta)$ is called $EY(\eta)$. For any $\eta' \in EY(\eta)$ we have $VY(\eta')=VY(\eta)$.

 \begin{lem} \label{lem:vraiPseudoPeriodique}
  Let $\Phi \in \Out(G)$ and let $f : T \to T$ be a representative for $\Phi$. Let $\eta$ be a pINP for $f$. Then for $g \in G$ the following are equivalent:
  \begin{enumerate}[label = (\roman*)] 
   \item $VY(\eta) = g VY(\eta)$
   \item $VY(\eta) \cap g VY(\eta) \neq \varnothing$
  \end{enumerate}

 \end{lem}
 \begin{proof}
  The implication (i) $\Rightarrow$ (ii) is obvious since $VY(\eta) \neq \varnothing$.
  
  Suppose (ii). Let us prove (i). Let $x \in VY(\eta) \cap g VY(\eta)$. Let $y \in VY(\eta)$. By definition of $VY(\eta)$, there is a concatenation of pINPs from $x$ to $y$. Since $x \in g VY(\eta)$ there is also a concatenation from $x$ to $gy$. Thus $y$ and $gy$ can be connected by a concatenation of pINP so they belong to the same Nielsen class. We deduce $g VY(\eta) = VY(\eta)$.
 \end{proof}

\begin{lem} \label{lem:constructionSousArbre}
 Let $\Phi \in \Out(G)$ with a train track representative $f : T \to T$. Let $\eta$ be a pINP for $\Phi$.
 
 Suppose $\stab(VY(\eta))$ contains a loxodromic element.
 
 The action of $\stab(VY(\eta))$ on the subtree $Y(\eta):= \conv(VY(\eta))$ is cocompact.
\end{lem}
\begin{proof} 
 Let $VY:= VY(\eta)$, $Y := \conv(VY)$ and $EY:= EY(\eta)$.
 
 The set $VY/\stab(VY)$ is finite because the image of $VY$ in $T/G$ is finite:  by Lemma \ref{lem:vraiPseudoPeriodique}, if $g \in G$,  $y \in VY$ and $gy \in VY$,  then $g \in \stab(VY)$. 
  
 The set $EY /\stab(VY)$ is finite because $EY / G$ is finite and for $g \in G$, $VY \cap gVY \neq \varnothing \Rightarrow g \in \stab(VY)$. There exist representatives $I_1, \dots, I_k$ for $EY$ and their union $I_1 \cup \dots \cup I_k$ is a compact subset of $T$ whose orbit covers $Y$.
 
 Thus $Y / \stab(VY)$ is compact.
\end{proof}

\begin{cor}
Let $\Phi \in \Out(G)$ with a train track representative $f : T \to T$ be such that there exists a pINP $\eta$ for $f$.
 The sets $VY(\eta)/\stab(VY(\eta))$ and $EY(\eta)/ \stab(VY(\eta))$ are finite.
\end{cor}

 The stabilizers of the Nielsen classes contain all the information on pseudo-periodic conjugacy classes, and they can be computed:
\begin{prop} \label{prop:calculSousArbres} 
 Let $\Phi \in \Out(G)$. Let $f : T \to T$ be a train track representative for $\Phi$.
 For every pINP $\eta$ for $f$, one can compute algorithmically a description of the set $VY(\eta)$ and of its stabilizer in the following form:
 \begin{itemize}
  \item a finite set of generators for the subgroup $\stab(VY(\eta))$ 
  \item a list of representatives for each $\stab(VY(\eta))$-orbit of points of $VY(\eta)$
  \item a list of representatives for each $\stab(VY(\eta))$-orbit of pINP in $EY(\eta)$.
 \end{itemize}
\end{prop}
\begin{proof} 
  By Lemma \ref{lem:compute-subdivision} we may subdivide $T$ such that the endpoints of all pINPs are vertices.
  
  Let $\eta$ be a pINP for $f$. 
  
  We construct the set $EY := EY(\eta)$ as follows. It suffices to construct all possible pINP concatenations starting from an endpoint of $\eta$, which actually constructs the vertices of $VY:=VY(\eta)$. 
  
 Define $EY_n$ as the list of pINPs which appear in concatenations of length at most $n$ from $y_0$, where $\eta = [y_0, y_1]$. For all $n \geq 1$ we have $EY_n \subset EY_{n+1}$.
 
 The construction of $EY_n$ is algorithmic. One needs a list of representatives for all $G$-orbits of pINPs, which can be computed using Proposition \ref{prop:compute-pinps}. Given a pINP $\gamma$ ending with a vertex $x$, we need to know whether there is a pINP $\gamma' \neq \gamma$ with $x$ as endpoint. This question can be answered using these facts, which are consequences of Proposition \ref{prop:algos-de-base}:
 \begin{itemize}
  \item given a $G$-orbit of pINPs $[\gamma']$ one can decide whether one endpoint of $[\gamma']$ is in the same orbit as $x$ 
  \item if so, one can compute a representative $\gamma'$ with $x$ as an endpoint and compute all other translates by applying repeatedly a generator of $G_x$. Since $T$ is locally finite, $G_x \cdot \gamma'$ is finite so eventually, all translates of $\gamma'$ starting at $x$ have been computed.
 \end{itemize}
 
 Since there are finitely many $G$-orbits of pINPs there exists $n \in \N$ such that all orbits of pINPs in $EY_{n+1}$ also appear in $EY_{n}$.

 Choose a minimal set of representatives in $\mathcal R \subset  EY_{n}$ such that all orbits of pINPs in $EY_{n+1}$ are represented in $\mathcal R$. Let $\eta_1, \dots, \eta_s$ be the pINPs of $EY_n \setminus \mathcal R$. 
 For every $1 \leq i \leq s$ there exists an element $g_i \in G$ such that $g_i \eta_i$ belongs to $\mathcal R \subset EY_{n-1}$.
Let $a$ be a generator of $\stab(\eta)$.
Let $G_\eta:= \langle g_1, \dots, g_s , a \rangle$.
 
 Then define $EY_\eta := G_\eta \cdot EY_n$. 
 
 \medskip
 
 Let us prove that  $EY_{\eta} = EY$. 
 It suffices to prove that $EY_{\eta}$ is the set of all pINPs which belong to a pINP path to $y_0$.  
 
 Suppose $\gamma \in EY_\eta$. There exists a word $g_1 \dots g_m$ in the generators of $G_\eta$ and a pINP $\gamma_0 \in EY_n$ such that $\gamma = g_1 \dots g_m \gamma_0$. 
 
 For every $1 \leq i \leq m$, the paths $\gamma_0, g_i \gamma_0$ belong to $EY_n$ and can be joined by a concatenation of pINPs within $EY_n$. Thus there is a concatenation containing all pINPs 
 \[
 \gamma_0, g_1 \gamma_0, g_1 g_2 \gamma_0, \dots, g_1 \dots g_m \gamma_0 = \gamma.
 \]
 
 Conversely suppose that there exists a concatenation of pINPs $\gamma_0, \gamma_1, \dots, \gamma_m$ such that $\gamma_0 \in EY_n$. Let us prove that $\gamma_m \in EY_\eta$. We proceed by induction on $m$. The case $m=0$ is obvious. For greater $m$ suppose every concatenation of at most $m$ pINPs starting with a pINP in $EY_n$ is contained in $EY_\eta$. By the induction hypothesis, there exists $g \in G_\eta$ such that  $g\gamma_{m-1}$ belongs to $EY_n$. Thus  $g\gamma_m$ belongs to $EY_{n+1}$ so there exists $i \in \{1, \dots, s\}$ such that $g_i g\gamma_m \in EY_n$. Since $g_i g \in G_\eta$, this proves that $\gamma_m \in EY_\eta$. Thus a concatenation of $m+1$ pINPs belongs to $EY_\eta$. We can conclude by induction.
 
 Now let us prove that the stabilizer of $EY_\eta$ is the subgroup $G_\eta$, which is finitely generated. The inclusion $G_\eta \subset \stab(EY_\eta)$ is true. Conversely let $g \in \stab(EY_\eta)$. There exists a concatenation of pINPs $\eta, \dots, g \cdot \eta$. There exists $g' \in G_\eta$ such that $g'\eta= g \cdot \eta$ so $g \in g'\stab(\eta) \subset G_\eta$.
 
 \bigskip
 
 The list $EY_n$ provides the list of representatives of the orbits of pINPs in $EY$. Moreover the ends of the elements of $EY_n$ are a list $VY_n$ of representatives of the orbits of $VY$. There may be redundant elements in both $EY_n$ and $VY_n$ and they may be eliminated algorithmically using Proposition \ref{prop:algos-de-base}.
\end{proof}

  Serre's Lemma (\cite[Section 6.5, Corollary 2]{SerreArbresAmalgames}) gives a criterion for deciding whether a finitely generated subgroup acting on a tree is elliptic.
\begin{lem} \label{lem:decideElliptique}
 Suppose a finitely generated group $G$ acts by isometries on a simplicial tree $T$. Let $\{s_1, \dots, s_n\}$ be a finite generating set for $G$. Then the action of $G$ is elliptic if and only if
 \begin{itemize}
  \item for all $i \in \{1, \dots, n\}$ the isometry $s_i$ is elliptic
  \item for all $i,j \in \{1, \dots, n\}$, $s_is_j$ is elliptic.
 \end{itemize}
\end{lem}
  
\begin{cor} \label{coro:computePeriodic}
 There is an algorithm which finds whether there exists a pseudo-periodic conjugacy class in $G$.
\end{cor}
\begin{proof}
 We can algorithmically compute finitely many generating sets for the stabilizers $\stab(VY(\eta))$ for finitely many pINPs $\eta$ representing every pINP orbit in $T$.
 
 By Lemma \ref{lem:vraiPseudoPeriodique} there exist pseudo-periodic conjugacy classes if and only if one of these subgroups contain a loxodromic element.
 
 Since a finitely generated subgroup is either elliptic or contains a loxodromic element, it then suffices to check that the algorithm returns only elliptic subgroups using Serre's lemma (Lemma \ref{lem:decideElliptique}). If it returns only elliptic subgroups then there exists no pseudo-periodic conjugacy class in $G$, otherwise there exists one.
\end{proof}

 \subsection{Link with pseudo-periodic subgroups}
 
  We now introduce \emph{pseudo-periodic subgroups}, which can be defined independently of any train track representative. We will prove  that when a train track map does exist, elements of pseudo-periodic conjugacy classes are pseudo-periodic for some $\psi$ in the outer automorphism class of $[\phi^k]$ for some $k \in \N$. Pseudo-periodic subgroups which contain loxodromic elements are contained in the subgroups $\stab (VY(\eta))$ defined before, and maximal pseudo-periodic subgroups containing loxodromic elements coincide with some $\stab(VY(\eta))$.

 \begin{defi}
  Let $\phi \in \Aut(G)$ be an actual automorphism. An element $g \in G$ is \emph{pseudo-periodic} for $\phi$ if for any $T \in \D$, for any $x \in T$, the sequence $\left (d_T(x, \phi^n(g)x) \right ) _{n \in \N}$ is bounded.
 \end{defi}
 \begin{remas} \label{rem:generales}
  \begin{enumerate}
   \item The boundedness of the sequence $\left (d_T(x, \phi^n(g)x) \right ) _{n \in \N}$ depends neither on $T$ nor on $x$.
   \item This notion is not defined for conjugacy classes of $g$. There may exist a pseudo-periodic element $g \in G$ and a conjugate of $g$ which is not pseudo-periodic. The conjugate will be pseudo-periodic for some $\phi' \in \Aut(G)$ such that $\phi'$ and $\phi$ belong to the same outer class.
   \item The definition also applies to elliptic elements. Not every elliptic element is pseudo-periodic in general. 
   \item If $g$ is pseudo-periodic and loxodromic, then its conjugacy class is pseudo-periodic for $\Phi=[\phi]$.
   \item The set $\{ g \in G/ g \text{ pseudo-periodic }\}$ is a subgroup of $G$. Indeed suppose $g,h$ are pseudo-periodic. For every $n \in \N$ we have 
   \begin{align*}
    d_T(x, \phi^n(gh)x) & \leq d_T(x, \phi^n(g)x) + d_T(\phi^n(g)x, \phi^n(g) \phi^n(h)x)\\
                        & = d_T(x, \phi^n(g)x) + d_T(x, \phi^n(h)x).
   \end{align*}
   \item A loxodromic element $g \in G$ is pseudo-periodic for $\phi$ if and only if there exists $m \in \N$ such that $g^m$ is pseudo-periodic. The direct implication is immediate. For the converse, observe that since $\phi^n(g)$ is loxodromic for every $n \in \N$ we have $d_T(x, \phi^n(g)x) \leq d_T(x, \phi^n(g^m)x)$. 
   \item  For any $n \in \N$ we have the inclusion $G_{\phi} \subset G_{\phi^n}$. 
  \end{enumerate} 
 \end{remas}
 
  \begin{defi}
  Let $\phi \in \Aut(G)$. The \emph{pseudo-periodic subgroup} $G_\phi$ associated to $\phi$ is the subgroup of pseudo-periodic elements for $\phi$.
 \end{defi}

 \begin{lem}\label{lem:pseudoPerHonnete}
  Let $\phi \in \Aut(G)$ and $g \in G_{\phi}$ be a loxodromic pseudo-periodic element. For each $T \in \D$, there exists $\in \in \N$ such that for every $k \in \Z$, $\phi^k(g)$ and $\phi^{k+n}(g)$ have a pair of fundamental domains in common.
 \end{lem}
 
 \begin{proof}
  Since $[g]$ is a pseudo-periodic conjugacy class, by Proposition \ref{prop:equivalencePseudoperiodique}, up to replacing $\phi$ by a power, we may suppose that $\|\phi^n(g)\|_T = \|g \|_T$ for every $n \in \N$.
  
  Let $x \in T$. Since $g$ is pseudo-periodic there exists $r > 0$ such that for all $n \in \N$, $d_T(x, \phi^n(g)x) \leq r$. Therefore, for every $n \in \N$, $x$ is at distance at most $r/2$ of the axis of $\phi^n(g)$. The intersection $\axe_T(\phi^n(g)) \cap \overline {B(x,2r)}$ is a segment with endpoints in $\overline{B(x, 2r)}$ and length at least $3r$.  Since there are finitely many such segments, there exist $n,m$ with $m > n + N$ such that $\axe_T(\phi^m(g)) \cap B(x,2r) = \axe_T(\phi^n(g)) \cap B(x,2r)$. Moreover we have $\|\phi^n(g)\|_T = \|g\|_T \leq r$. These axes overlap over more than $2 \|\phi^n(g)\|_T = 2 \|\phi^m(g)\|_T$. This means they share at least two consecutive fundamental domains in $T$.
  
  By Lemma \ref{lem:2domainesFonda}, $n$ does not depend on $T$: for any $S \in \D$, if $g, \phi^n(g)$ share two fundamental domains in $T$ then they also do in $S$. 
 \end{proof}

 \begin{rema} \label{rem:comportement-axes} 
  For a loxodromic pseudo-periodic element $g \in G$ we can predict two behaviours for $\phi^n(g)$ depending on the \emph{modulus} $\Delta(g) \in \Q^*$. The modulus is a morphism $\Delta: G \to \Q^*$ defined as follows; see \cite{LevittAutomorphisms} for a more detailed presentation. Fix an elliptic element $a \in G$. Since the commensurator of any elliptic element of $G$ is $G$, there exist $p, q \in \Z \setminus \{0\}$ such that $g a^p g^{-1} = a^q$. The ratio $p/q$ does not depend on the choice of $a$. Define $\Delta(g):= p/q \in \Q$.
  
  Suppose $g$ is pseudo-periodic. By Lemma \ref{lem:pseudoPerHonnete} there exists $n >1$ such that $\phi^n(g)$ and $g$ share a fundamental domain $\sigma$. Therefore $\phi^n(g)g^{-1}$ fixes one endpoint of $\sigma$. Let $a \in G$ be a generator of the stabilizer of this endpoint: then $\phi^n(g)g^{-1} \in \langle a \rangle$. There exists $k \in \Z$ such that $\phi^n(g) = a^k g$. 
  
  Let $p, q \in \Z \setminus \{0\}$ be such that $ga^pg^{-1} = a^q$.
  \begin{itemize}
   \item Suppose $\Delta(g) \neq 1$. Then $p-q \neq 0$ and $a^{p-q} g =a^p g a^{-p}$. Thus if $l_1 = l_2 \mod p-q$ then $a^{l_1}g$ and $a^{l_2} g$ are conjugate.
   
   Since there exist infinitely many possible choices for $n$, there exist $n, n'$ such that $\phi^n(g) = a^kg$ and $\phi^{n'} (g) = a^{k'}g$ with $k=k' \mod p-q$. So $\phi^n(g), \phi^{n'}(g)$ are conjugates, and so are $g$ and $\phi^{n'-n}(g)$. Thus the conjugacy class of $g$ is actually periodic. This does not imply that there are finitely many axes among the axes of $\{\phi^n(g) / n \in \N\}$!
   
   \item Suppose $\Delta(g) = 1$, so $q = p$. This is the \emph{unimodular} case. Then $a^p g = g a^p$ so $a^p$ fixes the axis of $g$. There exist $\phi^n(g) = a^kg, \phi^{n'}(g)= a^{k'}g$ with $k=k' \mod p$ so $\phi^n(g)$ and $\phi^{n'}(g)$ have the same axis. By Lemma \ref{lem:2domainesFonda} this implies that $g$ and $\phi^{n'-n}(g)$ have the same axis, and that $\phi^l(g)$ and $\phi^{l + n'-n} (g)$ have the same axis for every $l \in \Z$. In particular there exists $m >1$ such that the axis of $\phi^{mn} (g)$ is the same as the axis of $g$ for every $n \in \Z$. Here the axis is ``periodic'', but the conjugacy class of $g$ is not in general.
  \end{itemize}
 \end{rema}

 In the rest of the section, we assume that $\Phi \in \Out(G)$ is an automorphism of $G$ and it has a train track representative $f_0:T \to T$. We will study automorphisms in $\Aut(G)$ which are in the outer automorphism class $\Phi^k$ for some $k \in \N$. For $\psi$ in the outer automorphism class $\Phi^k$, we will always use the train track representative $f$ such that $f = g \cdot f_0^k$ for some $g \in G$, such that $f$ is $\psi$-equivariant. This way, periodic Nielsen paths are the same for every $\psi \in \Phi^k$.

 In our study, automorphisms whose representative has periodic points hold a special role.
\begin{defi}
 Let $\psi \in \Aut(G)$ be in the outer class of $\Phi^k$ for some $k \in \Z$ and $f:T \to T$ be a $\psi$-equivariant train track map. Let $\eta$ be a pINP for $f$. We say that $\psi$ is \emph{adapted to $\eta$} if the endpoints of $\eta$ are $f$-periodic.
\end{defi}

  \begin{rema} \label{rem:vy-puissance}
  The set of pINPs for $\Phi$ is the same as for $\Phi^n$, for any $n \in \N$. Thus a map $f$ is adapted to $\eta$ if and only if there exists $n \in \N$ such that $f^n$ is adapted to $\eta$. 
 \end{rema}

\begin{prop} \label{prop:implicationPseudo-pseudo}
 Let $k$ be such that for every pINP  $\eta$ in $T$, $[f_0^k(\eta)]$ is in the orbit of $\eta$.
 For every loxodromic $g \in G$ such that the conjugacy class $[g]$ is pseudo-periodic for $\Phi$, there exists $\psi \in \Phi^k$ with an associated train track representative $f$ such that $g \in G_{\psi}$ and $(\psi,f)$ is adapted to a pINP $\eta$ in the decomposition of the axis of $g$ given by Proposition \ref{prop:concatenationNielsen}.
 
 Therefore for every loxodromic $g$ in a pseudo-periodic conjugacy class $[g]$ of $G$, there exists a pINP $\eta$ and a pair $(\psi,f)$ adapted to $\eta$ such that $g \in G_{\psi}$. Moreover $g \in \stab(VY(\eta))$.
\end{prop}
\begin{proof}
 Let $[g]$ be a pseudo-periodic conjugacy class. By Proposition \ref{prop:concatenationNielsen} the axis of $g$ is a concatenation of pINPs. Take a pINP $\eta=[x, y]$ in $\axe_T(g)$: there exists $h\in G$ such that $[f_0^k(\eta)] = h \eta$. Define $\psi := c_{h^{-1}} \circ \phi^k$. The map $f:=h^{-1} \cdot f_0^k$ is a train track representative for $\psi$ and the endpoints $x, y$ of $\eta$ are fixed points for $f$.
 
 For every $n \in \N$, $f^n(x)=x$. Since $[x, gx]$ is a concatenation of pINPs, the distance $d_T(f^n(x), f^n(gx)) = d_T(x, \phi^n(g)x)$ is bounded so $g \in G_\psi$.
\end{proof}

 A point $x \in T$ is a \emph{pre-periodic point} for $f$ if there exists $n \geq 0$ such that $f^n(x)$ is periodic. We say $x$ is \emph{non-escaping} if for any $y \in T$ the distance $d_T(y, f^n(x))$ is bounded for $n \rightarrow \infty$. A pre-periodic point is obviously non-escaping but the converse is also true for a train track map: 
  \begin{lem} \label{lem:preperImpliquePer}
  Let $\phi \in \Aut(G)$ and $f:T \to T$ be a train track representative for $\phi$. If $x$ is non-escaping then there exists $n \in \N$ and $m \in \N \setminus \{0\}$ such that $f^{n+m}(x) = f^n(x)$, so $x$ is actually a pre-periodic point.
 \end{lem}
 \begin{proof}
  Suppose for some $n_0 \in \N$, $f^{n_0}(x)$ is a vertex. Since $f$ maps vertices to vertices and $T$ is locally finite, $f^n(x)$ takes only finitely many values so there exists $n< m$ such that $f^n(x) = f^m(x)$ and $f^n(x)$ is an actual periodic point.
  
  If $f^n(x)$ is in the interior of an edge for every $n \in \N$, there exists $n <m$ such that $f^n(x)$ and $f^m(x)$ are both contained in the same edge $e$. The image $f^{m-n}(e)$ contains $e$ so there exists a fixed point $x_0 \in e$ for $f^{m-n}$. Since $e$ is a legal segment, it is stretched uniformly by a factor $\lambda > 1$ and unless $f^n(x)= x_0$, the distance $d_T(x_0, f^{k(m-n)}(x))$ grows exponentially when $k \rightarrow + \infty$. Therefore $f^n(x)=x_0$.
 \end{proof} 
 
 Suppose that the pair $(\phi, f)$ is adapted to a pINP $\eta \subset T$.
 The points of $VY(\eta)$ are non-escaping so they are also pre-periodic for $f$.

 The following lemma is a general result about actions on trees.
 \begin{lem} \label{lem:dihedral} 
 Let $H$ be a group acting on a simplicial tree $T$ by isometries such that there exist loxodromic isometries. Suppose that there exists an elliptic element $a \in H$ such that for every loxodromic $h \in H$, the product $ah$ is an elliptic isometry. Then the action of $H$ on $T$ is dihedral.
\end{lem}
\begin{proof}
 As in \cite{GuirardelLevitt07} we distinguish 4 possible types of action for $H$: linear abelian, dihedral, genuine abelian and irreducible. We prove the lemma by contraposition, i.e. if the action is abelian or irreducible then for every elliptic $a$ there exists a loxodromic $h \in H$ such that $ah$ is loxodromic.
 
 Suppose that the action is abelian (linear or genuine). It has a fixed point $\xi \in \partial T$. Let $a \in H$ be an elliptic isometry. Let $h \in H$ be a loxodromic isometry. Both $h$ and $a$ fix $\xi$ so there exists $x \in \axe_T(h)$ such that the subray $[x, \xi]$ is contained in $\Fix_T(a)$. The isometry $a h$ acts on $[x, \xi]$ like $h$ so it is also loxodromic.
 
 Otherwise suppose that the action is irreducible. Let $a \in H$ be an elliptic element, and $h \in H$ be a loxodromic element.
 Suppose $\axe_T(h) \cap \Fix_T(a) = \varnothing$. Let us prove that $ah$ is loxodromic. Let $p \in T$ be the projection of $\Fix_T(a)$ on $\axe_T(h)$. 
 The point $p' :=h^{-1}p$ is sent by $ah$ to the point $ap$. Let $v \in T$ be a point in $[p',p]$. Then $ahv \in [ap, ahp]$.  In particular $ahv$ is not in $[p', ahp']$ so $ah$ is loxodromic and $p'$ belongs to its axis.
 
 Suppose now that $\axe_T(h) \cap \Fix_T(a)$ is non empty. If $a$ fixes one of the endpoints of $\axe_T(h)$ then the same discussion as for the abelian case applies. Else there exists a loxodromic element $g \in H$ whose axis does not intersect the axis of $h$. Since $\axe_T(h)$ is not fixed by $a$, there exists $k \in \Z$ such that $h^k \axe_T(g)$ does not intersect $\Fix_T(a)$ so up to replacing $g$ by the conjugate $h^k g h^{-k}$ we may assume $\axe_T(g) \cap \Fix_T(a) = \varnothing$. Then $ag$ is loxodromic.
\end{proof}

Suppose that a pair $(\phi, f)$ is adapted to $eta$.
The following is the key for proving that $VY(\eta)$ is stable by $G_\phi$:
\begin{lem} \label{lem:connecterAx}
 Suppose the group $G_\phi$ of pseudo-periodic elements contains a loxodromic element.
 Let $\phi \in \Aut(G)$ with a train track representative $f: T\to T$. Suppose that $(\phi, f)$ is adapted to $\eta = [x,x']$ a pINP in $T$. Then $G_\phi \subset \stab(VY(\eta))$.
\end{lem}

\begin{figure}
  \centering
  \includegraphics[scale=0.6]{\localhost/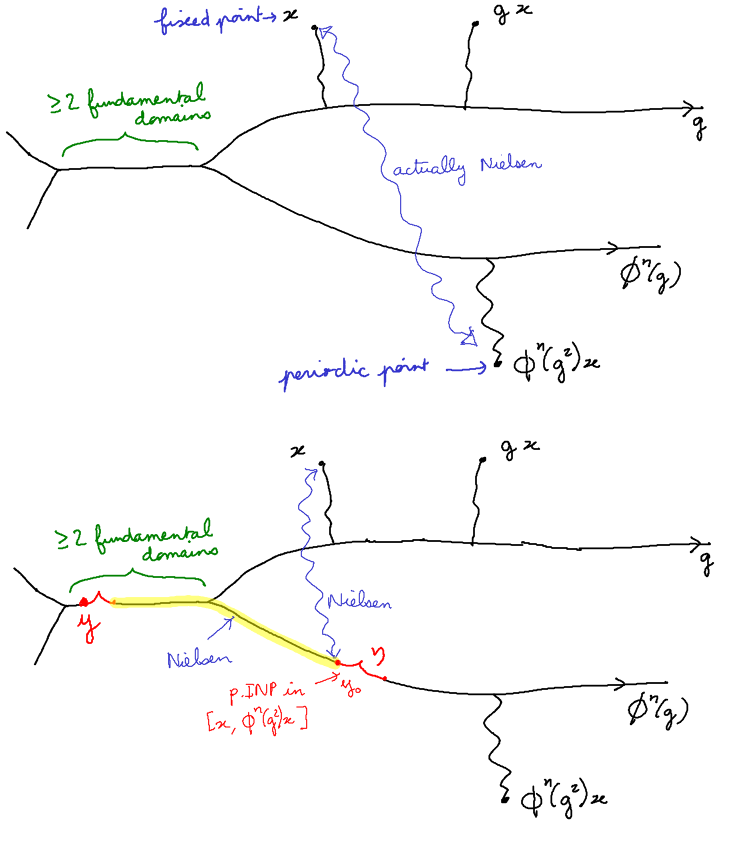}
  \caption{Finding a concatenation of pINPs between $x$ and $gx$} \label{fig:raboutageINP}
 \end{figure}

\begin{proof}
 Suppose first that $g$ is a loxodromic element.
 By Remark \ref{rem:vy-puissance}, $f$ is adapted to $\eta$ if and only if its powers are. By Remark \ref{rem:generales} $G_{\phi} \subset G_{\phi^n}$ so it suffices to prove the lemma for the pair $(\phi^n, f^n)$ for some $n \in \N$. Thus we may assume for simplicity that $f(x)=x$.
 
 By Lemma \ref{lem:vraiPseudoPeriodique} the points $gx$ and $g^2x$ are pre-periodic. They might not be periodic points. If $gx$ is a periodic point, the path $[x, gx]$ is a periodic Nielsen path so by Lemma \ref{lem:NielsenInclus} it is a concatenation of pINP.

 Suppose otherwise. The axis of $g$ is a legal concatenation of pINPs, by Proposition \ref{prop:concatenationNielsen}. We will prove that all these pINPs belong to $EY(\eta)$. Then by Lemma \ref{lem:vraiPseudoPeriodique} $g \in \stab(VY(\eta))$.
 This will conclude the proof for $g$ loxodromic.

 See Figure \ref{fig:raboutageINP} for a picture of this proof.
 
 The point $x$ is fixed by $f$ and since the path between $x$ and $g^2x$ is a concatenation of pINPs, the point $g^2x$ is non-escaping, so it is pre-periodic. Thus there exists $N \in \N$ such that for all $n \geq N$, $f^n(g^2x) = \phi^n(g^2)x$ is a periodic point.  By Lemma \ref{lem:pseudoPerHonnete} we may choose $n$ such that $\phi^n(g^2)x$ is periodic and the axes of $\phi^n(g)$ and $g$ intersect along two common fundamental domains. By Lemma \ref{lem:vraiPseudoPeriodique} we have $\phi^n(g^2) \in \stab(VY(\eta))$.
 
 The path $[x, \phi^n(g^2)x]$ is a concatenation of pINPs. It contains a pINP $\eta' \subset \axe_T(\phi^n(g))$. By Lemma \ref{lem:NielsenInclus} $\eta'$ appears in the decomposition of the Nielsen path $[x, \phi^n(g^2)x]$. Thus $\eta' \in EY(\eta)$.
 
 Since $\phi^n(g)$ and $g$ share two fundamental domains, there exists $l \in \Z$ such that $\phi^n(g)^l \eta' \subset \axe_T(g) \cap \axe_T(\phi^n(g))$. By Lemma \ref{lem:NielsenInclus} this translate of $\eta'$ appears in the pINP decomposition of $\axe_T(g)$, besides it belongs to $EY(\eta)$. 
 This completes the proof when $g$ is loxodromic.
 
 \bigskip
 
 If $g$ is elliptic the argument using the axes does not work any more. If there exist loxodromic elements $h_1, h_2$ such that $g = h_1 h_2$ then the path $[x, gx]$ can be written as the concatenation $[x, h_1x] \cdot h_1[x, h_2x]$ and the loxodromic case completes the proof.
 
 If $g$ cannot be written as a product of loxodromic elements then by Lemma \ref{lem:dihedral} the action of $G_\phi$ on its minimal subtree is dihedral (by assumption it is not elliptic). Let $\ell$ be the invariant axis in $T$. If $g$ fixed $\ell$ then for any loxodromic element $h \in G_\phi$, $gh$ would be loxodromic so $g=(gh)h^{-1}$ would be a product of loxodromic elements. Thus $g$ acts like a symmetry. There exists a loxodromic isometry $h \in G_\phi$ with axis $\ell$. By Proposition \ref{prop:concatenationNielsen} $\ell$ is a legal concatenation of pINPs. 
 
 As in the loxodromic case, the pINPs in the decomposition of $\ell= \axe_T(h)$ belong to $EY(\eta)$.
 Furthermore this decomposition into pINPs is $g$-invariant. Let $y \in VY(\eta)$ be an endpoint of a pINP in $\ell$, 
 then the point $gy$ is also in $EY(f)$. By Lemma \ref{lem:vraiPseudoPeriodique} we have $g \in \stab(VY(\eta))$. 
 This concludes the proof in the elliptic case.
\end{proof}

 \begin{cor} \label{coro:gphi=stabvy}
  Let $\phi \in \Aut(G)$  and let $f$ be a train track representative for $\phi$. Suppose that 
  $(\phi, f)$ is adapted to $\eta$.  
  Then $\stab(VY(\eta)) = G_{\phi}$.
 \end{cor}
 \begin{proof}
  The inclusion $G_{\phi} \subset \stab(VY(\eta))$ is proved by Lemma \ref{lem:connecterAx}. Let us prove the inverse inclusion. Suppose $g \in \stab(VY(\eta))$. Let $x_0$ be an endpoint of $\eta$. Then $gx_0 \in VY(\eta)$ so $gx_0$ is pre-periodic. For all $n \in \N$ we have $f^n(g x_0) = \phi^n(g) f^n(x_0)$. Since $x_0$ and $gx_0$ are pre-periodic, there exists $B > 0$  such that the sets $\{f^n(x_0), n \in \N\}$ and $\{\phi^n(g) f^n(x_0), n \in \N\}$ have diameter smaller than $B$. Then for all $n \in \N$ we have 
  \[
   d_T(x_0, \phi^n(g) x_0) \leq d_T(x_0, \phi^n(g) f^n(x_0)) + \phi^n(g) d_T(f^n(x_0), x_0) \leq 2B
  \]
  so $g \in G_\phi$.
 \end{proof}

Pseudo-periodic subgroups defined by pairs $(\phi, f)$ which are adapted to a pINP have good properties but not all pseudo-periodic subgroups arise in this way. However we will see that it is the case for \emph{maximal} pseudo-periodic subgroups.

\begin{prop} \label{prop:action-irreductible}
  Let $\psi \in [\phi^k]$ for some $k \in \N$, and let $f$ be the associated train track representative. If the action of $G_\psi$ on $T$ is irreducible then the pair $(\psi, f)$ is adapted to a pINP.
\end{prop}
\begin{proof}
 Suppose that the action of $G_{\psi}$ is irreducible. Then we claim that $G_{\psi}$ contains unimodular loxodromic elements $u,v \in G$ whose axes have distinct endpoints and an intersection longer than $4C_f$. We postpone the proof of the claim.
 
 By Remark \ref{rem:comportement-axes} there exists $n \geq 1$ such that $u$ and $\phi^n(u)$ have the same axis with same orientation, and have the same action on it.
 Recall that $f^n$ sends maximal legal segments to legal segments. 
 Writing the axis of $u$ as a bi-infinite concatenation of pINPs $\{\eta_i, i \in \Z\}$, $f^n$ induces a translation on the set of pINPs: there exists $k_u \in \Z$ such that $[f^n(\eta_i)] = \eta_{i + k_u}$.
 
 Similarly $f^n$ shifts the pINPs $\{\eta'_j, j \in \Z\}$ in $\axe_T(v)$ by $k_v \in \Z$. 
 The length of a pINP is at most $2C_f$ so the intersection $\axe_T(v) \cap \axe_T(u)$ contains a pINP $\eta_i= \eta'_j$ belonging to the decomposition of both axes. Suppose by contradiction that $k_u \neq 0$ is non-zero, then there exists $m \in \N$ such that $ [f^{nm}(\eta_{i})]\cap \axe_T(u) = \eta_{i+mk_u}$, which is non-empty, does not intersect $\axe_T(v)$. However $[f^{nm}(\eta_i)] = f^{nm}(\eta'_j) = \eta'_{j + mk_v} \subset \axe_T(v)$. This is a contradiction.
 
 Thus $f^n(\eta_i)= \eta_i$ so $(\psi,f)$ is adapted to $\eta_i$.
 
 \bigskip
 
 Now let us prove the claim. We will use Lemma \ref{lem:elements-axes-disjoints} below to construct $u$ and $v$ as announced. Since the action of $G_\psi$ is irreducible there exist loxodromic elements $g,h \in G_\psi$ with disjoint axes. The axes of $g$ and $hg^{-1}h^{-1}$ are disjoint so the element $u:=ghg^{-1} h^{-1}$ is unimodular and loxodromic, and $\|u\|_T \geq 2 \|g\|_T$. Up to replacing $g$ by a power we may suppose $\|u\|_T \geq 4 C_f$. Moreover the axes of $g$ and $u$ intersect along one fundamental domain of $g$, so there exists $l \in \Z$ such that the axes of $u$ an $g^lu^{-1} g^{-l}$ do not intersect. Thus the element $v:= g^l u^{-1} g^{-l}u$ is also unimodular and loxodromic, moreover the axes of $u$ and $v$ overlap on a length equal to $\|u\|_T$, which proves the claim.
\end{proof}

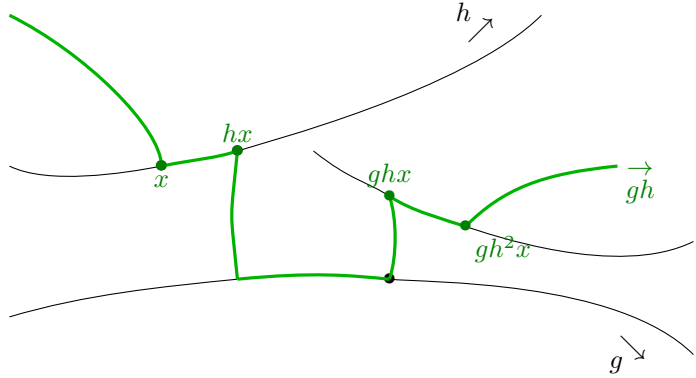
\begin{figure}
 \centering
  %\usetikzlibrary{calc}
\begin{tikzpicture}

 \newcommand{\bez}[4]{ %début, fin, vecteur 1, vecteur 2
    \draw #1 .. controls ($ #1 + #3$) and ($#2 + #4$) .. #2;
}
 \newcommand{\bezv}[4]{ %début, fin, vecteur 1, vecteur 2
    \draw[green!70!black, very thick] #1 .. controls ($ #1 + #3$) and ($#2 + #4$) .. #2;
}

%\bez{(0,0)}{(3,0)}{(0,1)}{(0,1)};

\coordinate (ag) at (0,0);
\coordinate (bg) at (9,-0.5);
\coordinate (ah) at (0,2);
\coordinate (bh) at (7,4);
\coordinate (x) at (2,2);
\coordinate (y) at (3,0.5);
\coordinate (hx) at (3,2.2) ;

\coordinate (gy) at (5,0.5);

\draw (x) node {$\bullet$};
\draw (hx) node {$\bullet$};
\draw (gy) node {$\bullet$};

\bez{(ag)}{(y)}{(1,0.3)}{(-1,-0.1)};
\bezv{(gy)}{(y)}{(-0.5,0.05)}{(1,0.1)};

\bez{(bg)}{(gy)}{(-1,1)}{(0.5,-0.05)};
%\bez{(ah)}{(bh)}{(2,-0.5)}{(-1,-1)};

\bez{(ah)}{(x)}{(0.5,-0.25)}{(-0.5,-0.1)};
\bezv{(hx)}{(x)}{(-0.3,-0.1)}{(0.5,0.1)};
\bez{(hx)}{(bh)}{(0.3,0.1)}{(-1,-1)};

\bezv{(hx)}{(y)}{(-0.1,-1)}{(-0.1,1)};

%axe ghg
\coordinate (gah) at (4,2.2);
\coordinate (gbh) at (9,1);
\coordinate (ghx) at (5,1.6);
\coordinate (gx) at (6,1.2) ;

\bez{(gah)}{(ghx)}{(0.5,-0.4)}{(-0.3,0.2)};
\bezv{(gx)}{(ghx)}{(-0.3,0.1)}{(0.3,-0.2)};
\bez{(gx)}{(gbh)}{(0.3,-0.1)}{(-1,-0.5)};

\bezv{(ghx)}{(gy)}{(0.1,-0.3)}{(0.1,0.3)};

\bezv{(x)}{(0,4)}{(0,0.5)}{(1,-0.5)}
\bezv{(gx)}{(8,2)}{(0.5,0.5)}{(-1,-0.1)}

\draw [green!50!black] (8.3,1.95) node {$\rightarrow$} node [below] {$gh$};
\draw (6.2,3.8) node {$\nearrow$} node [above left] {$h$};
\draw (8.2,-0.4) node {$\searrow$} node [below left] {$g$};

\draw[green!50!black] (x) node {$\bullet$} node [below] {$x$};
\draw[green!50!black] (hx) node {$\bullet$} node [above] {$hx$};
\draw[green!50!black] (ghx) node {$\bullet$} node [above] {$ghx$};
\draw[green!50!black] (gx) node {$\bullet$} node [below right] {$gh^2x$};

\end{tikzpicture}
 \caption{Axis of $gh$ when $g,h$ have disjoint axes} \label{fig:axes-disjoints}
\end{figure}

\begin{lem}\label{lem:elements-axes-disjoints}
 Let $T$ be a simplicial tree. Let $a,b$ be loxodromic isometries of $T$ with disjoint axes. Then the product $ab$ is a loxodromic isometry with translation length $\|a\|_T + \|b\|_T + 2 d_T(\axe_T(a), \axe_T(b))$. Moreover $\axe_T(ab) \cap \axe_T(a)$ is a single fundamental domain for $a$.
\end{lem}
\begin{proof}
 Most of this is proved in \cite[1.5]{CullerMorganGroupActions}; we give the proof for completeness but it is best explained by Figure \ref{fig:axes-disjoints}.
 Let $a, b$ be as above. Let $p_a$ be the projection of $\axe_T(b)$ to $\axe_T(a)$, and $p_b$ be the projection of $\axe_T(a)$ to $\axe_T(b)$. Let $x:= b^{-1}p_b$. Then $abx = a p_b$ and the path $[x, abx]$ is the union $[x, p_b] \cup [p_b, p_a] \cup [p_a, a p_a] \cup [ap_a, ap_b]$. The overlap of those 4 segments has zero length because $p_a, p_b$ are projections. The path $[p_b, ab p_b] = [p_b, p_a] \cup [p_a, ap_a] \cup [ap_a, ap_b] \cup [ap_b, ab p_b]$ has the same length since $[ap_b, abp_b] = ab[x, p_b]$ and the decomposition also has trivial overlap. If $x$ were not in the axis, the point $p_b$ which lies in the interior of $[x, abx]$ would be moved by a shorter distance, so $x$ is actually in the axis. The translation length is $\|a\|_T + \|b\|_T + 2 d_T(p_a, p_b)$ as announced and $\axe_T(ab) \cap \axe_T(a) = [p_a, a p_a]$.
\end{proof}

Pseudo-periodic subgroups containing loxodromic elements are contained in stabilizers of Nielsen classes, with equality when the subgroup comes from a pair adapted to a pINP:
\begin{lem}\label{lem:pas-nice-contenu-nice}
 For every $\psi$ in $\Phi^k$ for $k \in \N$, if $G_{\psi}$ contains a loxodromic element, there exists $\eta \subset T$ a pINP such that $G_{\psi} \subset \stab(VY(\eta))$.
 
 If $(\psi,f)$ is adapted to a pINP then the inclusion is an equality.
\end{lem}
\begin{proof}
 The second statement, for $(\psi,f)$ adapted to a pINP $\eta$, is a consequence of Corollary \ref{coro:gphi=stabvy}: we have  $G_{\psi}  = \stab(VY(\eta))$.

 Suppose $G_{\psi}$ contains a loxodromic element. 
 If the action of $G_{\psi}$ on $T$ is irreducible then by Proposition \ref{prop:action-irreductible} $(\psi,f)$ is adapted to a pINP $\eta$ so the second statement applies.
 
 Otherwise suppose that there is a fixed point in $\partial T$ for the action of $G_{\psi}$. Let $g$ be such that the attracting point for $g$ is fixed by $G_{\psi}$. 
 
 By Proposition \ref{prop:implicationPseudo-pseudo} there exists a pINP $\eta \subset \axe_T(g)$ and a pair $(\psi', f')$ adapted to $\eta$ so $g \in \stab(VY(\eta)) =  G_{\psi'}$. Let $h \in G_{\psi}$ be a loxodromic element. The axes of $g$ and $h$ intersect along an infinite ray. There exist powers $l,m$ such that $\|g^l\|_T = \|h^m\|_T$ so $g^l$ and $h^m$ have at least two fundamental domains in common, and by Lemma \ref{lem:2domainesFonda} for every $\phi$ in the outer class of $\Phi$ and every $n \in \Z$ the elements $\phi^n(g^l)$ and $\phi^n(h^m)$ also share two fundamental domains.

 Up to replacing $\eta$ by some translate $g^j \eta$, we may assume $\eta \subset \axe_T(h)$ so $g^lh^{-m}\eta = \eta$. Since the conjugacy class $[h]$ is pseudo-periodic, $\axe_T(h)$ is a concatenation of pINPs and by Lemma \ref{lem:NielsenInclus}, $\eta$ is a pINP of the decomposition. Thus the pINPs in $\axe_T(h)$ belong to $EY(\eta)$. By Lemma \ref{lem:vraiPseudoPeriodique} $h \in \stab(VY(\eta))$.
 
 By Lemma \ref{lem:dihedral} all elliptic elements can be written as the product of loxodromic elements so since all loxodromic elements of $G_\psi$ belong to $\stab(VY(\eta))$, we have $G_{\psi} \subset \stab(VY(\eta))$.
 
 \bigskip

The only remaining case is when the action of $G_{\psi}$ is dihedral.
In that case $G_{\psi}$ preserves an axis $\ell$ in $T$. 

There exists a loxodromic element $g \in G_{\psi}$ whose axis is $\ell$. Let $\eta$ be a pINP in $\ell$: then the endpoints of the pINPs in the decomposition of $\ell$ belong to $VY(\eta)$. Since this decomposition is unique by Lemma \ref{lem:NielsenInclus} and $G$ preserves the set of pINPs in $T$, the group $G_{\psi}$ preserves $VY(\eta) \cap \ell$. By Lemma \ref{lem:vraiPseudoPeriodique} this implies $G_{\phi} \subset \stab(VY(\eta))$. 
\end{proof}

 Thus maximal pseudo-periodic subgroups which contain loxodromic elements come from pairs adapted to a pINP:
\begin{prop} \label{prop:nice-implique-maximal}
 A pseudo-periodic subgroup $H$ containing a loxodromic element is maximal for inclusion among pseudo-periodic subgroups if and only if there exists a pINP $\eta \subset T$ and a pair $(\psi, f)$ adapted to $\eta$ such that $H = G_{\psi}$.
\end{prop}
\begin{proof}
 Suppose $H$ is maximal. By Lemma \ref{lem:pas-nice-contenu-nice} there exists a pINP $\eta$ such that $H \subset \stab(VY(\eta))$. There exists $\psi \in \Phi^k$ for some $k \in \N$ and $g \in G$ such that $f:= gf_0^k$ fixes the endpoints of $\eta$, so the pair $(\psi, f)$ is adapted to $\eta$ and $\stab(VY(\eta)) = G_{\psi}$. By maximality $H = G_{\psi}$.
 
 Conversely suppose there exists a pair $(\psi, f)$ such that $H = G_{\psi}$ adapted to a pINP $\eta \subset T$. By Corollary \ref{coro:gphi=stabvy} $H = \stab(VY(\eta))$. Let $H'$ be a pseudo-periodic subgroup such that $H \subset H'$: by Lemma \ref{lem:pas-nice-contenu-nice} there exists $\eta'$ such that $H' \subset \stab(VY(\eta'))$. Thus for $h \in H$, $\axe_T(h) \subset \conv(VY(\eta'))$. By uniqueness of the decomposition of $\axe_T(h)$ into pINPs, the pINPs in $h$ belong to $EY(\eta)$ and $EY(\eta')$ so $VY(\eta)= VY(\eta')$ and $H' \subset \stab(VY(\eta)) =H$, hence the maximality of $H$.
\end{proof}

\begin{cor} \label{cor:nbFiniClassesConjugaisonSousgroupes}
 There are finitely many conjugacy classes of maximal pseudo-periodic subgroups containing loxodromic elements associated to the outer classes $\Phi^k \in \Out(G)$ for $k \in \N$. Furthermore, these subgroups are finitely generated.
\end{cor}
\begin{proof}
This is a consequence of Proposition \ref{prop:nice-implique-maximal} and Corollary \ref{coro:gphi=stabvy}: the maximal pseudo-periodic subgroups containing loxodromic elements coincide with the stabilizers of Nielsen classes containing loxodromic elements.

The conclusion follows from the fact that there are finitely many orbits of Nielsen classes and their stabilizers are finitely generated.
\end{proof}

\section{Algorithm} \label{sec:conclu}

In this section we assemble the results obtained before to give the desired algorithm. First let us explain how to deal with restricted deformation spaces.

Let $G$ be a GBS group represented by a graph of groups $\Gamma$. As in \cite{PapinWhitehead} we say that a family of subgroups $\A$ which is invariant by conjugacy and taking subgroups is \emph{represented} by finite sets of integers $(I_v)_{v \in V(\Gamma)}$ if for every $v \in \Gamma$ and every lift $\tilde v$ in the universal cover of $\Gamma$, the set $I_v$ is the set of minimal elements of $\{[G_v : G_v \cap A], A \in \A\}$ for divisibility.

\begin{theo}
 There is an algorithm which takes
 \begin{itemize}
  \item a non-elementary GBS group $G$
  \item an automorphism $\phi \in \Aut(G)$
  \item a train track representative $f:T \to T$ for $\phi$ where $T$ is given as a graph of groups
 \end{itemize}
 which finds out whether $\phi$ is pseudo-atoroidal.
\end{theo}
\begin{proof}
  This is a consequence of Corollary \ref{coro:computePeriodic}.
\end{proof}

\begin{theo}
 There is an algorithm which takes
 \begin{itemize}
  \item a non-elementary GBS group $G$
  \item a pseudo-atoroidal automorphism $\phi \in \Aut(G)$
  \item a train track representative $f:T \to T$ for $\phi$ where $T$ is given as a graph of groups
  \item a family of sets $(I_v)_{v \in V(T/G)}$ representing the family of subgroups $\A$
 \end{itemize}
 which finds out whether $\phi$ has no pseudo-periodic element. In that case it decides whether $\phi$ is iwip for $\D^\A$.
\end{theo}
In the previous theorem, $\phi$ need not be pseudo-atoroidal as long as it has no simple pseudo-periodic conjugacy class. Note that we do not know how to check this condition algorithmically.
\begin{proof}
 There are two steps in this algorithm.
 
 The first step consists in finding either a primitive train track representative or a proof of reducibility for $\phi$. Corollary \ref{coro:construireTTirred} solves this problem algorithmically using the matrix $A(f)$ and produces a primitive train track representative if it exists. 
 
 The second step consists in computing turns of the stable lamination and its Whitehead graphs, then using Corollary \ref{coro:feuille}. Using Lemma \ref{lem:algoFeuilles} we can compute all Whitehead graphs. Then compute connected components of $\Wh_T(\Lambda^+,v)$ for every $v \in V(T/G)$. Let $W$ be such a Whitehead graph. If $W$ is not connected we need to compute the stabilizers of connected components.
 As in \cite[Section 3]{PapinWhitehead}, one can compute for every connected component $C$ the index $i(C):= [G_v : \stab(C)]$. Then $\stab(C) \in \A$ if and only if $i(C)$ is divisible by some $i \in I_v$.

 If there exists a Whitehead graph with a connected component whose stabilizer is in $\A$ then by Proposition \ref{prop:WhNonConnexeImpliqueRed} $\phi$ is reducible. Moreover the proof of the proposition also gives an invariant class of special factors with respect to $\A$ for $\phi$.

 Otherwise Corollary \ref{coro:feuille} states that $\phi$ is fully irreducible. 
\end{proof}

  \bibliography{../bibli_iwips}
 \end{document}